\documentclass[11pt]{article}
\usepackage{latexsym}
\usepackage{amsfonts,amssymb,amsmath}
\newtheorem{thm}{Theorem}[section]
\newtheorem{cor}[thm]{Corollary}
\newtheorem{lem}[thm]{Lemma}
\newtheorem{pro}[thm]{Proposition}
\newtheorem{defn}[thm]{Definition}

\title{Bernoulli measure on strings, and Thompson-Higman monoids }

\author{ Jean-Camille Birget  }

\date{\today}
\begin{document}
\maketitle

\begin{abstract}
The Bernoulli measure on strings is used to define {\em height functions} 
for the dense ${\cal R}$- and ${\cal L}$-orders of the Thompson-Higman 
monoids $M_{k,1}$.
The measure can also be used to characterize the ${\cal D}$-relation
of certain submonoids of $M_{k,1}$.
The computational complexity of computing the Bernoulli measure of certain 
sets, and in particular, of computing the ${\cal R}$- and ${\cal L}$-height
of an element of $M_{k,1}$ is investigated.
\end{abstract}

%%%%%%%%%%%%%%%%%%%%%%%%%%%%%%%%%%%%%%%%%%%%%%%%%%%%%%%%
% Section 1
%%%%%%%%%%%%%%%%%%%%%%%%%%%%%%%%%%%%%%%%%%%%%%%%%%%%%%%%

\section{Introduction}

Since their introduction in the 1960s \cite{Th0, McKTh, Th} the groups of
Richard J.\ Thompson have become well known for their remarkable 
properties. They were generalized by Graham Higman \cite{Hig74} to a class
of groups $G_{k,i}$ ($k \geq 2$, $k > i \geq 1$) with similar properties: 
They are finitely presented infinite simple groups, containing all finite
groups. The Thompson-Higman groups have a large literature; see the list of 
references in \cite{BiThomps, BiThomMon}. 
The groups $G_{k,i}$ can be generalized in a straightforward way to monoids
$M_{k,i}$ which also have remarkable properties:
They are finitely generated with word problem decidable in polynomial time, 
they are congruence-simple, and they contain all finite monoids 
\cite{BiThomMon}, \cite{BiRL}, \cite{BiJD}. 

The Green relations of $M_{k,1}$  were characterized in \cite{BiRL},
\cite{BiJD}: The ${\cal J}$- and the ${\cal D}$-relations do not have much 
structure, as $M_{k,1}$ is ${\cal J}$-0-simple and has exactly $k-1$ 
non-zero ${\cal D}$-classes. 
The ${\cal R}$- and ${\cal L}$-orders are complicated:
They are dense, and are characterized by image sets and partitions, in 
analogy with the monoid of all partial functions on an infinite set. We 
will define $M_{k,1}$ and describe its Green relations in more detail below. 
Here we will only talk about $M_{k,1}$ ($k \geq 2$), and ignore 
$M_{k,i}$ for $i \neq 1$. 

The aim of this paper is to introduce {\it height functions} for the 
${\cal R}$- and ${\cal L}$-orders of $M_{k,1}$, which is not obvious for 
dense orders.
It turns out that the Bernoulli measure on the free monoid $A^*$ (over a 
finite alphabet $A$ with $|A| = k$) plays a crucial role for defining 
${\cal R}$- and ${\cal L}$-height functions in $M_{k,1}$.
The Bernoulli measure enables us also to characterize the 
${\cal D}$-relation of certain submonoids of $M_{k,1}$, namely submonoids 
of transformations that map fixed-length words to fixed-length words. 
We also determine the computational complexity of computing the Bernoulli
measure associated with an element of $M_{k,1}$.
The appendix contains a proof that $M_{k,1}$ is congruence-simple (the 
proof of this fact in \cite{BiThomMon} was incomplete).

%%%%%%%%

\subsection{Definition of the Thompson-Higman groups and monoids}

In order to make this paper (mostly) self-contained we start with 
definitions and notations. We will define $M_{k,1}$ by a partial
action on $A^*$, the free monoid over $A$ (consisting of all finite 
sequences of elements of $A$), where $A$ is an alphabet of size $k$. 
We call the elements of $A^*$ {\it words}, and we also include the 
{\it empty word} $\varepsilon$. 
We denote the {\it length} of $w \in A^*$ by $|w|$.
The {\it concatenation} of $u,v \in A^*$ is denoted by $uv$ or by
$u \cdot v$; more generally, the concatenation of $B, C \subseteq A^*$ is
$BC = \{uv : u \in B, v \in C\}$.
A word $u \in A^*$ is a {\it prefix} of $v \in A^*$ iff $uz = v$ for
some $z \in A^*$, and we write $u \ {\sf pref} \ v$.  We say that $u$ and
$v$ are {\it prefix-comparable} if $u \ {\sf pref} \ v$ or
$v \ {\sf pref} \ u$.
A {\it prefix code} is a subset $C \subseteq A^*$ whose elements are 
two-by-two prefix-incomparable.
A prefix code is maximal iff it is not strictly contained in any other 
prefix code. 

A subset $R \subseteq A^*$ is called a {\it right ideal} iff 
$RA^* \subseteq R$. We call $R$ an {\it essential} right ideal iff $R$
intersects every right ideal of $A^*$. More generally, for right ideals $R'
\subseteq R \subseteq A^*$, $R'$ is {\em essential in} $R$ iff $R'$ 
intersects all right ideals included in $R$. 

A right ideal $R$ is {\it generated} by a set $C \subseteq A^*$ iff $R$ 
is the intersection of all right ideals that contain $C$; equivalently, 
$R = CA^*$.
One proves easily that a right ideal $R$ has a unique minimal (under
inclusion) generating set, and that this minimal generating set is a 
{\it prefix code}. 
The prefix code that generates $R$ is maximal iff $R$ is an essential 
right ideal.
In this paper we only use right ideals that are {\it finitely generated}.

For a partial function $f:A^* \to A^*$ we denote the domain by
${\sf Dom}(f)$ and the image by ${\sf Im}(f)$. 
A function $\varphi: R_1 \to A^*$ is a {\em right ideal homomorphism} of 
$A^*$ iff ${\sf Dom}(\varphi) = R_1$ is a right ideal and for all 
$x_1 \in R_1$ and all $w \in A^*$: \ $\varphi(x_1w) = \varphi(x_1) \ w$. 
It follows that ${\sf Im}(\varphi)$ is a right ideal, and if $R_1$ is 
finitely generated (as a right ideal) then ${\sf Im}(\varphi)$ is finitely 
generated.
We write the action of partial functions on the left of the argument;
equivalently, functions are composed from right to left.

A right ideal homomorphism $\varphi: R_1 \to R_2$ is uniquely determined
by its restriction $P_1 \to S_2$, where $P_1$ is the prefix code
that generates $R_1$ as a right ideal, and $S_2$ is a set (not necessarily 
a prefix code) that generates $R_2$ as a right ideal. 
This finite total surjective function $P_1 \to S_2$ is called the 
{\it table of} $\varphi$. The finite prefix code $P_1$ is called the
{\em domain code} of $\varphi$ and is denoted by ${\sf domC}(\varphi)$. 
When $S_2$ is a prefix code it will be denoted by ${\sf imC}(\varphi)$ and
called {\em image code}.

A right ideal homomorphism $\Phi: R'_1 \to A^*$ is called an 
{\em essentially equal restriction} of a right ideal homomorphism
$\varphi: R_1 \to A^*$ (or, equivalently, $\varphi$ is an {\em essentially
equal extension} of $\Phi$) iff $R'_1$ is essential in
$R_1$, and for all $x'_1 \in R'_1$: \ $\varphi(x'_1) = \Phi(x'_1)$.
The multiplication in $M_{k,1}$ (and in $G_{k,1}$) depends on the 
following facts: \ (1) Every homomorphism $\varphi$ between finitely 
generated right ideals of $A^*$ has a {\em unique} maximal essentially equal 
extension (Prop.\ 1.2(2) in \cite{BiThomMon}); this extension is denoted 
by ${\sf max}(\varphi)$. 
 \ (2) Every right ideal homomorphism $\varphi$ has an essentially equal 
restriction $\varphi'$ whose table $P' \to Q'$ is such 
that both $P'$ and $Q'$ are prefix codes (remark after Prop.\ 1.2 in 
\cite{BiThomMon}).
 
We are now ready to define the {\it Higman-Thompson monoid} $M_{k,1}$: 
As a set, $M_{k,1}$ consists of all homomorphisms (between finitely 
generated right ideals of $A^*$) that have been maximally essentially 
equally extended. In other words, as a set,

\smallskip

 \ \ \ $M_{k,1} = \{ {\sf max}(\varphi) : \varphi$ is a homomorphism
between finitely generated right ideals of $A^* \}$.

\smallskip

\noindent
The multiplication is composition followed by maximal essentially equal
extension (which is unique). This multiplication is associative (Prop.\ 1.4 
in \cite{BiThomMon}). Thus we have a {\em partial action} of $M_{k,1}$ on 
$A^*$ (since the multiplication is not just composition, but needs to be 
followed by maximal essential extension).

The Higman-Thompson monoid $M_{k,1}$ also has a true action by partial
functions on $A^{\omega}$ (the Cantor space, consisting of all
$\omega$-sequences over $A$). The action of $\varphi \in M_{k,1}$ on 
$z \in A^{\omega}$ is defined by $\varphi(z) = y w$ if $z$ can be written
as $z = xw$ for some $x \in {\sf domC}(\varphi)$, where $y = \varphi(x)$; 
$\varphi(z)$ is undefined if $z$ has no prefix in ${\sf domC}(\varphi)$. 

For a right ideal $R \subseteq A^*$ generated by a prefix code $P$
we call $P A^{\omega}$ the set of {\em ends} of $R$, denoted by
${\sf ends}(R)$. We call two right ideals $R_1, R_2$ {\em essentially equal} 
iff ${\sf ends}(R_1) = {\sf ends}(R_2)$, and we denote this by 
$R_1 =_{\sf ess} R_2$.

%%%
\subsection{ The Bernoulli measure} 

For a fixed alphabet $A$ with $|A| = k$, and any set $X \subseteq A^*$  
we define the {\em Bernoulli measure} of $X$ by
$$\mu(X) \ = \ \sum_{x \in X} k^{-|x|} \, .$$
For the empty set the measure is 0. When $X$ is infinite, $\mu(X)$ can be 
any positive real number or $+ \infty$. 
For a non-empty finite set $X$, $\mu(X)$ it is a strictly positive $k$-ary 
rational number. 

By definition, a {\it $k$-ary rational number} is a number of the form 
$a / k^n$, where $a \in {\mathbb Z}$ and $n \in {\mathbb N}$; equivalently, 
it is a rational number that has a finite representation in base $k$. The 
ring of $k$-ary rational numbers is denoted by $\mathbb{Z}[\frac{1}{k}]$. 
For any $r = a / k^n \in \mathbb{Z}[\frac{1}{k}]$ we say that $a / k^n$ 
is {\it $k$-reduced} iff $a$ is not divisible by $k$; 
in that case we denote the numerator $a$ by ${\sf num}(r)$.

Obviously, the definition of $\mu$ amounts to viewing each word in $A^*$ 
as a sequence of independent $k$-ary Bernoulli trials in which each
choice $a_i \in A$ has the same probability, namely $\frac{1}{k}$.

If $X$ is a prefix code (finite or infinite) then \ $\mu(X) \leq 1$, with 
equality iff $X$ is maximal. This is the {\it Kraft inequality} (or 
equality). Hence $\varphi \in M_{k,1}$ is {\it total} (on the Cantor space 
$A^{\omega}$) iff $\mu({\sf domC}(\varphi)) = 1$;
$\varphi$ is {\it surjective} (onto $A^{\omega}$) iff 
$\mu({\sf imC}(\varphi)) = 1$.

\begin{lem} \label{set_rewrite} \   
For any prefix code $P \subset A^*$, \ $x \in P$, and $n \geq 0$, we have:
 \ \ $\mu(P) = \mu\big( (P - \{x\}) \ \cup \ x A^n\big).$
%More generally, for any maximal prefix code $Q \subset A^*$,
% \ $\mu(P) = \mu\big( (P - \{x\}) \cup x Q\big)$.
\end{lem}
{\bf Proof.} This is easy to verify, based on the fact that 
 \ $k^{-|x|} \ = \ \sum_{w \in A^n} \ k^{-|xw|}$ . 
  \ \ \ $\Box$

\begin{lem} \label{classReplacement} (from \cite{BiRL}). 
 \ Let $R_1 = P_1A^*$ and $R_2 = P_2A^*$ be right ideals, where $P_1$ and
$P_2$ are finite prefix codes. Then $R_1 =_{\sf ess} R_2$ iff $P_2$
can be transformed into $P_1$ by a finite sequence of replacements steps
of the following form:

\smallskip

\noindent {\bf (r1)} \ \ \ For a finite prefix code $C$ and for $c \in C$,
replace $C$ by \ $(C - \{c\}) \cup cA$.

\smallskip

\noindent {\bf (r2)} \ \ \ For a finite prefix code $C'$ such that
$cA \subseteq C'$ for some word $c$, replace $C'$ by \ $(C' - cA) \cup \{c\}$.
\end{lem}
{\bf Proof.} This is straightforward. \ \ \ $\Box$

\begin{pro} \label{ess_equal_same_mu} \   
If $P_1, P_2 \subset A^*$ are any prefix codes and
$P_1A^* =_{\sf ess} P_2A^*$, then $\mu(P_1) = \mu(P_2)$.
\end{pro}
{\bf Proof.} This follows easily from Lemmas \ref{set_rewrite} and
\ref{classReplacement}  above.           \ \ \ $\Box$

%%%%%%%%%%%%%%%%%%%%%%%%%%%%%%%%%%%%%%%%%%%%%%%%%%%
%% Section
%%%%%%%%%%%%%%%%%%%%%%%%%%%%%%%%%%%%%%%%%%%%%%%%%%%

\section{Height functions for the $\cal R$- and $\cal L$-orders 
         of $M_{k,1}$}

In a finite monoid, the ${\cal R}$-height of an element $x$ is defined to
be the length of a longest ascending strict ${\cal R}$-chain from a minimal
${\cal R}$-class to $x$ (and similarly for the ${\cal L}$-height).
For infinite monoids there may be no good way to define an ${\cal R}$- or 
${\cal L}$-height function at all, especially if the ${\cal R}$- and 
${\cal L}$-orders are dense (as is the case for $M_{k,1}$, by Section 4 in 
\cite{BiRL}). The main result of this section is that for $M_{k,1}$ we can
nevertheless define an ${\cal R}$-height and an ${\cal L}$-height. 
For a pre-order, in general, we define height functions as follows:

\begin{defn} \label{Heightfunct}
 \ A {\it height function} for a pre-order $(M, \preccurlyeq)$ is any
function $h: M \to {\mathbb R}$ such that for all $x, y \in M$:
 \ if $y \preccurlyeq x$ then $h(y) \leq h(x)$, and 
if $y \prec x$ then $h(y) < h(x)$ \ (where, $y \prec x$ means 
``$y \preccurlyeq x$ and $x \not\preccurlyeq y$'').
It follows that $h(x) = h(y)$ when $x \equiv y$ (where $x \equiv y$ means
``$x \preccurlyeq y$ and $y \preccurlyeq x$'').
\end{defn}

The definition of the ${\cal R}$- and ${\cal L}$-height functions will be
guided by the characterizations of the ${\cal R}$- and ${\cal L}$-orders in
$M_{k,1}$, given in \cite{BiRL}.
Regarding $\leq_{\cal R}$ we have for all $\psi, \varphi \in M_{k,1} :$ 

\smallskip

$\bullet$ \ \ $\psi \ \leq_{\cal R} \ \varphi$ \ \ iff

\smallskip

$\bullet$
 \ \ ${\sf ends}({\sf Im}(\psi)) \subseteq {\sf ends}({\sf Im}(\varphi))$
 \ \ iff

\smallskip

$\bullet$ \ \ every right ideal of $A^*$ that intersects
${\sf Im}(\psi)$ also intersects ${\sf Im}(\varphi)$.

\smallskip

The characterization of $\leq_{\cal L}$ is more complicated and involves 
right-congruences. We first need some definitions. 
For any right ideal homomorphism $\varphi$, let ${\sf part}(\varphi)$ be
the right congruence on ${\sf Dom}(\varphi)$ defined by
$(x_1,x_2) \in {\sf part}(\varphi)$ iff $\varphi(x_1) = \varphi(x_2)$.
This definition of ${\sf part}(\varphi)$ can be extended to
${\sf ends}({\sf Dom}(\varphi))$: for $v_1, v_2 \in A^{\omega}$ we have 
$(v_1,v_2) \in {\sf part}(\varphi)$ iff $\varphi(v_1) = \varphi(v_2)$; this 
is iff there exist $w \in A^{\omega}$ and $x_1,x_2 \in {\sf Dom}(\varphi)$ 
such that
$(x_1,x_2) \in {\sf part}(\varphi)$ and $v_1 = x_1 w$, $v_2 = x_2 w$.

Two right congruences $\simeq_1$ (on $P_1A^*$) and $\simeq_2$ (on $P_2A^*$) 
are called {\it essentially equal} iff 
${\sf ends}(P_1A^*) = {\sf ends}(P_2A^*)$ and the extension of 
$\simeq_1$ to ${\sf ends}(P_1A^*)$ is equal to the extension of $\simeq_2$ 
to ${\sf ends}(P_1A^*)$. We denote this by 
$\simeq_1 \ =_{\sf ess} \ \simeq_2$. We say that $\simeq_1$ is an 
{\it essentially equal extension} of $\simeq_2$ (and equivalently, 
$\simeq_2$ is an {\it essentially equal restriction} of $\simeq_1$) iff 
$P_2A^* \subseteq P_1A^*$ and $\simeq_1 \ =_{\sf ess} \ \simeq_2$.

Let $\simeq_{{\sf domC}(\varphi)}$ be the restriction of the right 
congruence ${\sf part}(\varphi)$ to the finite prefix code 
${\sf domC}(\varphi)$.
We would like ${\sf part}(\varphi)$ to be determined in a simple way by 
$\simeq_{{\sf domC}(\varphi)}$ as follows.
Let $P \subset A^*$ be a finite prefix code and let $\simeq_P$ be an
equivalence relation on $P$; we call a right congruence $\simeq$ on 
$PA^*$ a {\it prefix code congruence} (determined by $\simeq_P$)
iff for all $p_1, p_2 \in P :$

\smallskip

$\bullet$ \ if $p_1 \simeq_P p_2$ then for all $w \in A^*$,
 \ $p_1 w \simeq p_2 w$ ; 

\smallskip

$\bullet$ \ for all $x,y \in A^*:$ \ if \ $p_1 \not\simeq_P p_2$ \ or if
 \ $x \neq y$ \ then \ $p_1 x \not\simeq p_2 y$ .

\smallskip

\noindent
In other words, $\simeq$ is a prefix code congruence on $PA^*$ iff $\simeq$
is the smallest (i.e., finest) right congruence that agrees with the 
restriction of $\simeq$ to $P$.

Although ${\sf part}(\varphi)$ is always a right congruence, it is not 
always a prefix congruence. In \cite{BiRL} it was proved that 
${\sf part}(\varphi)$ is a prefix congruence (determined by 
$\simeq_{{\sf domC}(\varphi)}$) \ iff \ $\varphi({\sf domC}(\varphi))$ 
is a prefix code.
In that case $\varphi({\sf domC}(\varphi))$ is denoted by 
${\sf imC}(\varphi)$. 

Let $\simeq$ be a prefix code congruence on $PA^*$ determined by 
$\simeq_P$. 
If $C \subseteq P$ is a class of $\simeq_P$ then a 
{\it class-wise replacement step} consists of replacing $C$ by the set of 
classes $\{C a_1, \ldots, C a_k\}$ and replacing $P$ by 
 \ $Q = (P - C) \, \cup \, C a_1 \, \cup \ \ldots \ \cup \, C a_k$ , where 
$A = \{a_1, \ldots, a_k\}$.
In \cite{BiRL} it was proved that the resulting equivalence relation 
$\simeq_Q$ also determines a prefix code congruence (on $QA^*$) which is 
essentially equal to the congruence determined by $\simeq_P$.
Hence, if we apply a finite sequence of class-wise replacement steps or 
inverses of replacement steps to a prefix code congruence $\simeq$, we 
obtain a prefix code congruence that is essentially equal to $\simeq$.
Conversely, it was proved in \cite{BiRL} that if two prefix code 
congruences $\simeq_1$ and $\simeq_2$ are essentially equal then each one
is obtained from the other one by a finite number of class-wise replacement
steps and their inverses. 

A prefix code congruence $\simeq$ is called {\it maximal} iff $\simeq$
is maximal with respect to $\subseteq_{\sf ess}$. 
It is easy to see that inverse class-wise replacements form a terminating 
and confluent rewriting system. Hence, every prefix code congruence $\simeq$ 
is $\subseteq_{\sf ess}$-contained in a unique maximal prefix code 
congruence, which we denote by ${\sf max}(\simeq)$. 
 
\smallskip

We can now state the characterization of the $\cal L$-order of
$M_{k,1}$ given in \cite{BiRL}.  For $\varphi, \psi \in M_{k,1}$,

\smallskip

$\bullet$ \ \ $\psi \ \leq_{\cal L} \ \varphi$ \ \ iff

\smallskip

$\bullet$ \ \ every class of ${\sf ends}({\sf part}(\psi))$
is a union of classes of ${\sf ends}({\sf part}(\varphi))$ \ \ iff

\smallskip

$\bullet$ \ \ every class of ${\sf part}(\psi)$ is a union of classes of
${\sf max}({\sf part}(\varphi))$.

%\smallskip

%%%
\subsection{A height function for the $\cal R$-order }

We just saw that the ${\cal R}$-order in $M_{k,1}$ is determined by the 
inclusion relation between the sets ${\sf ends}({\sf Im}(\varphi))$.
This suggests the following definition of an ${\cal R}$-height function 
for $M_{k,1}$. 

\begin{defn} \label{def_Rheight} \
Let $\varphi \in M_{k,1}$ be described by a table $P \to Q$ where
$P, Q \subset A^*$ are finite prefix codes, i.e., $P = {\sf domC}(\varphi)$
and $Q = {\sf imC}(\varphi)$.
Then the {\bf $\cal R$-height} of $\varphi$, denoted by 
${\sf height}_{\cal R}(\varphi)$, is defined by
$${\sf height}_{\cal R}(\varphi) \ = \ \mu({\sf imC}(\varphi)).$$
\end{defn}
We will prove below (Prop.\ \ref{Rheight_vs_Rorder}) that 
${\sf height}_{\cal R}$ depends only on the $\cal R$-class of $\varphi$ and
that it is indeed a height function for $(M_{k,1}, \, \leq_{\cal R})$. 
In particular, it does not depend on the particular table $P \to Q$ used 
to represent $\varphi$.

We said that we only use tables $P \to Q$ where $Q$ is a prefix code. Let 
us briefly investigate what happens when $Q$ is not a prefix code.

\begin{pro} \label{muQ_notprefcode} \   
Let $\Phi \in M_{k,1}$ be described by a table $P \to Q$, where $Q$ is 
{\em not} a prefix code, and let $\varphi$ be any essentially equal 
extension or restriction of $\Phi$ such that ${\sf imC}(\varphi)$ is a 
prefix code. Then, \ $\mu(Q) > \mu({\sf imC}(\varphi))$. 
\end{pro}
{\bf Proof.} \ When $Q$ is not a prefix code then there exists a prefix 
code $Q_0 \subset Q$ such that $QA^* = Q_0A^*$. In fact,
$Q_0 = Q - Q A A^*$, so $Q_0$ is uniquely determined by $Q$.  
Since $Q_0A^* = QA^* =_{\sf ess} {\sf imC}(\varphi) \, A^*$, and since 
$Q_0$ and ${\sf imC}(\varphi)$ are prefix codes, we have (by Prop.\ 
\ref{ess_equal_same_mu}): \ $\mu(Q_0) = \mu({\sf imC}(\varphi))$.
Moreover, $\mu(Q) = \mu(Q_0) + \mu(Q - Q_0) > \mu(Q_0)$. 
 \ \ \ $\Box$

\bigskip

\noindent The following example illustrates Prop.\ \ref{muQ_notprefcode}.

\medskip

\noindent {\sf Example.} \ Let $A = \{a,b\}$ and let $\varphi \in M_{2,1}$ 
be given by the following tables, all describing the same element of 
$M_{2,1}$:

\medskip

$\Phi_1 = $  \  \begin{tabular}{|c|c|c|}
$aa$  & $ab$   & $b$ \\ \hline
$a$   & $aa$   & $aaa$ \
\end{tabular} \ , \ \ \ $\Phi_2 = $ \ 
\begin{tabular}{|c|c|c|c|}
$aaa$ & $aab$ & $ab$   & $b$ \\ \hline
$aa$  & $ab$  & $aa$   & $aaa$ \
\end{tabular} \ , 

\medskip

$\Phi_3 = $ \ 
\begin{tabular}{|c|c|c|c|c|}
$aaaa$ & $aaab$ & $aab$ & $ab$   & $b$ \\ \hline
$aaa$  & $aab$  & $ab$  & $aa$   & $aaa$ \
\end{tabular} \ , \ and
 \ \ $\Phi_4 = $ \ 
\begin{tabular}{|c|c|c|c|c|c|}
$aaaa$ & $aaab$ & $aab$ & $aba$ & $abb$ & $b$  \\ \hline
$aaa$  & $aab$  & $ab$  & $aaa$ & $aab$ & $aaa$ \! \
\end{tabular} \ .   

\medskip

\noindent The measure of the set $\{a,aa,aaa\}$ in the table of
$\Phi_1$ is \ $\mu(\{a,aa,aaa\}) = \frac{7}{8}$. 
By essentially equal restriction we obtain $\Phi_2$ 
and the measure of the image set of the table of $\Phi_2$ is 
 \ $\mu(\{aa,ab, aaa\}) = \frac{5}{8}$.
Next we get the table $\Phi_3$ with measure
 \ $\mu(\{aaa,aab,ab,aa\}) = \frac{3}{4}$.
Finally, another restriction step yields $\Phi_4$, for which
$\Phi_4({\sf domC}(\Phi_4)) = {\sf imC}(\Phi_4))$ is a prefix code. Thus
we finally obtain the measure 
 \ $\mu({\sf imC}(\Phi_4)) = \mu(\{aaa,aab,ab, aaa, aab, aaa\})$ 
$ = \mu(\{aaa,aab,ab\})$ $ = \frac{1}{2}$; 
 \ so, $\mu({\sf imC}(\varphi)) = \frac{1}{2}$.
%Observe also that ${\sf Im}(\Phi_1)$ has $\{a\}$ as a minimum generating
%set, and \ $\mu(\{a\}) = \frac{1}{2}$. 
 \ \ \ $\Box$

\bigskip

\noindent As a consequence we have:

\begin{pro} \label{mu_characteriz_prefcongr} \   
For a right-ideal homomorphism $\varphi$ the following are equivalent:

\smallskip

\noindent {\bf (1)} \ ${\sf part}(\varphi)$ is a prefix congruence;

\smallskip

\noindent {\bf (2)} \ $\varphi({\sf domC}(\varphi))$ is a prefix code;

\smallskip

\noindent {\bf (3)} \     
$\mu(\varphi({\sf domC}(\varphi))) \ = \ $
${\sf min}\{ \mu(Q) : \ P \to Q$ is a table that represents $\varphi$\}.
\end{pro}
{\bf Proof.} The equivalence of (1) and (2) was proved in \cite{BiRL}.
The equivalence of (2) and (3) is given by Prop.\ \ref{muQ_notprefcode},
combined with Prop.\ \ref{ess_equal_same_mu}.
  \ \ \ $\Box$

\medskip

\noindent The term ``$\cal R$-height'' is justified by the following. 

\begin{pro} \label{Rheight_vs_Rorder} \
For all $\varphi \in M_{k,1}$, \ $\mu({\sf imC}(\varphi))$ depends
only on the $\cal R$-class of $\varphi$.

The function ${\sf height}_{\cal R}: \varphi \in M_{k,1} \ \longmapsto $
$ \ \mu({\sf imC}(\varphi)) \in {\mathbb Z}[\frac{1}{k}] \cap [0,1]$ 
 \ is a height function for the pre-order \ $(M_{k,1}, \leq_{\cal R})$ 
 \ (according to Definition \ref{Heightfunct}). 
\end{pro}
{\bf Proof.} \ The fact that $\mu({\sf imC}(\varphi))$ depends only on 
the $\cal R$-class of $\varphi$ follows immediately from the 
characterization of the $\cal R$-order mentioned at the beginning of 
section 2  (Theorem 2.1 in \cite{BiRL}), and from Prop.\ 
\ref{ess_equal_same_mu} above.

Suppose $\psi$ and $\varphi$ are represented by tables
$P_{\psi} \to Q_{\psi}$, respectively $P_{\varphi} \to Q_{\varphi}$.
We only prove that $\varphi >_{\cal R} \psi$ iff
${\sf height}_{\cal R}(\varphi) > {\sf height}_{\cal R}(\psi)$ (under the
assumption that $\psi$ and $\varphi$ are $\cal R$-comparable); the other
two case have a similar proof.
By Theorem 2.1 in \cite{BiRL}, $\varphi >_{\cal R} \psi$ iff
${\sf ends}(Q_{\psi}A^*) \subsetneqq {\sf ends}(Q_{\varphi}A^*)$. After an
essentially equal restriction of $\psi$ and $\varphi$, if necessary,
the latter holds iff $Q_{\psi} \subsetneqq Q_{\varphi}$.
This holds iff $\mu(Q_{\psi}) < \mu(Q_{\varphi})$, under the assumption
that $Q_{\psi}$ and $Q_{\varphi}$ are comparable under inclusion.
 \ \ \ $\Box$

\medskip

The following Lemma and Proposition show that the height function 
${\sf height}_{\cal R}$ is {\it onto} $\mathbb{Z}[\frac{1}{k}] \cap [0,1]$,
and that for any chain of numbers in $\mathbb{Z}[\frac{1}{k}] \cap [0,1]$
there are corresponding $<_{\cal R}$-chains.
\begin{lem} \label{prefcodesformeansure} \!\!. \\   
{\bf (1)} \ For every $h \in \mathbb{Z}[\frac{1}{k}] \cap [0,1]$ there 
exists a finite prefix code $P_h \subset A^*$ with $\mu(P_h) = h$.

\smallskip

\noindent {\bf (2)}  \ For all $g, h \in \mathbb{Z}[\frac{1}{k}]$ 
with \ $0 \leq g < h \leq 1$ \ the finite prefix codes $P_g, P_h$ 
constructed in {\bf (1)} satisfy 

\smallskip

 \ \ \ \ \  $P_g \subset P_h A^*$, and 
 \ $P_g A^* \neq_{\sf ess} P_h A^*$. 

\smallskip

\noindent Hence the set 
 \ $\{P_h A^* : h \in \mathbb{Z}[\frac{1}{k}] \cap [0,1] \}$ \ is a chain 
of right ideals, no two of which are essentially equal, and with 
$\mu(P_h) = h$. 
\end{lem} 
{\bf Proof.} 
{\bf (1)} \ As before, $A = \{a_1, \ldots, a_k\}$.
When $h = 0$ we pick $P_0 = \varnothing$, and when $h = 1$ we pick 
$P_1 = \{ \varepsilon\}$ (where $\varepsilon$ is the empty string). 
We assume next that $0 < h < 1$. Then $h$ has a unique finite base-$k$
expansion of the form \ $h = 0{\bf .}d_1 \ldots d_i \ldots d_n$ \ where
$d_n \neq 0$, and \ $d_i \in \{0, 1, \ldots, k-1\}$ \ for $i = 1, \ldots, n$.
With $h$ we associate the following prefix code:

\medskip

$P_h \ = \ \{a_1, \ldots, a_{d_1}\}$   $ \ \ \cup \ \ $
$ a_{d_1 +1} \, \{a_1, \ldots, a_{d_2}\}$   $ \ \ \cup \ \ $
$ a_{d_1 +1} a_{d_2 +1} \, \{a_1, \ldots, a_{d_3}\}$
$ \ \cup \ \ $  $ \ \ldots \ $

\smallskip

\hspace{.4in}  $ \ldots $  $ \ \ \cup \ \ $  
$a_{d_1 +1} \ldots a_{d_{i-2}+1} \, \{a_1, \ldots, a_{d_{i-1}}\}$
$ \ \ \cup \ \ $
$a_{d_1+1} \ldots a_{d_{i-2}+1} a_{d_{i-1}+1} \, \{a_1, \ldots, a_{d_i}\}$

\smallskip

\hspace{.4in} $ \ \ \cup \ \ $ 
 $a_{d_1 +1} \ldots a_{d_{i-1}+1} a_{d_i+1} \, \{a_1, \ldots,a_{d_{i+1}}\}$ 
 $ \ \ \cup \ \ $  $ \ldots $

\smallskip

\hspace{.4in}
$ \ldots $   $ \ \ \cup \ \ $
$a_{d_1+1} \ldots a_{d_{n-1}+1} \, \{a_1, \ldots, a_{d_n}\}$ .

\medskip

\noindent Here, the set $\{a_1, \ldots, a_{d_i}\}$ is empty if $d_i = 0$
($1 \leq i \leq n$).
It is easy to verify that 

\smallskip

 \ $\mu(a_{d_1 +1} \ldots a_{d_{i-1} +1} \{a_1, \ldots, a_{d_i} \}) \ = $
$ \ k^{-i} \, d_i$,  

\smallskip

\noindent hence $h = \mu(P_h)$. 
Note also that \ $|P_h| = \sum_{i=1}^n d_i$. 

\medskip

\noindent For example, for $k = 5$ and $h = 0.0031042$ (in base 5)
we have

\smallskip

$P_h \ = \ a_1 a_1 \, \{a_1,a_2,a_3\} \ \cup \ $
$a_1 a_1 a_4 \, \{a_1\} \ \cup \  a_1 a_1 a_4 a_2 a_1 \, \{a_1,a_2,a_3,a_4\}$
$ \ \cup \ a_1 a_1 a_4 a_2 a_1 a_5 \, \{a_1,a_2\}$. 

\medskip

\noindent {\bf (2)}  \ The above construction of $P_h$ from $h$ has the 
following property:
For any $k$-ary rationals $h, g$ with $0 < h < g \leq 1$ we have
$P_h \subset P_gA^*$. This is proved next.

When $g =1$ this is clear (since then $P_g$ is a maximal prefix code, and 
$P_h$ is not maximal).
When $h < g < 1$ then $g$ has a base-$k$ expansion of the form
 \ $g = 0.d_1 \ldots d_{i-1} b_i \ldots b_m$ with $d_i < b_i$
 \ (for some $i$ with $1 \leq i \leq {\sf min}\{m,n\}$), where 
$b_i, \ldots, b_m$ are base-$k$ digits. Then

\medskip

$P_g \ = \ \{a_1, \ldots, a_{d_1}\}$   $ \ \ \cup \ \ $
$ a_{d_1 +1} \, \{a_1, \ldots, a_{d_2}\}$   $ \ \ \cup \ \ $
$ a_{d_1 +1} a_{d_2 +1} \, \{a_1, \ldots, a_{d_3}\}$
$ \ \cup \ \ $  $ \ \ldots \ $

\smallskip

\hspace{.4in}  $ \ \ldots \ $  $ \ \ \cup \ \ $
$a_{d_1 +1} \ldots a_{d_{i-2}+1} \, \{a_1, \ldots, a_{d_{i-1}}\}$
$ \ \ \cup \ \ $
$a_{d_1 +1} \ldots a_{d_{i-2}+1} a_{d_{i-1}+1} \, \{a_1,$
$\ldots,a_{d_i},\ldots,a_{b_i}\}$

\smallskip

\hspace{.4in} $ \ \ \cup \ \ $
$a_{d_1+1}\ldots a_{d_{i-1}+1} a_{b_i+1}\, \{a_1,\ldots, a_{b_{i+1}}\}$
  $ \ \cup \ \ $  $ \ \ldots \ \ \ $

\hspace{.4in} 
$ \ldots \ \ $   $ \ \cup \ $
$a_{d_1 +1} \ldots a_{b_{m-1}+1} \, \{a_1,\ldots, a_{b_m}\}$ .

\medskip

\noindent Since $d_i < b_i$ we have

\medskip

$a_{d_1 +1} \ldots a_{d_{i-1}+1} \, \{a_1, \ldots, a_{d_i}\}$
 \ \ $\subset$ \ \   
$a_{d_1+1}\ldots a_{d_{i-1}+1}\, \{a_1,\ldots,a_{d_{i}},\ldots,a_{b_{i}}\}$.

\medskip

\noindent And we have for all $j$ ($i \leq j \leq n$):

\medskip

$a_{d_1 +1} \ldots a_{d_{i-1}+1} \ \ \ \ \ a_{d_i+1} \ \ \ \ \ a_{d_{i+1}+1}$
$\ldots $
$a_{d_{j-1}+1} \, \{a_1, \ldots, a_{d_j} \}$ $ \ \ \subset \ \ $ 

$a_{d_1 +1} \ldots a_{d_{i-1}+1}\, \{a_1,\ldots, a_{b_{i}}\}\, A^*$;

\medskip

\noindent indeed, $a_{d_i +1} \in \{a_1,\ldots, a_{b_{i}}\}$ (because 
$d_i +1 \leq b_i$), and 
 \ $a_{d_{i+1}+1} \ldots a_{d_{j-1}+1} \, \{a_1, \ldots, a_{d_j} \}$ 
$ \subset A^*$.
Thus, $P_h \subset P_gA^*$, hence
$P_hA^* \subset P_gA^*$. Also, $P_hA^* \neq_{\sf ess} P_gA^*$ since
$\mu(P_hA^*) = h \neq g = \mu(P_gA^*)$.
 \ \ \ $\Box$ 

\medskip

\noindent {\sf Notation:}  For a subset $X \subseteq A^*$, ${\sf id}_X$
denotes the partial identity function. In other words, for any $w \in A^*$, 
${\sf id}_X(w) = w$ if $w \in X$; and ${\sf id}_X(w)$ is undefined if
$w \not\in X$.

\begin{pro} \label{every-rat-is_Rheight} \ {\bf (1)} 
 \ For every $h \in {\mathbb Z}[\frac{1}{k}] \cap [0,1]$ there exists 
$\varphi_h \in M_{k,1}$ such that ${\sf height}_{\cal R}(\varphi_h) = h$.
Moreover, $\varphi$ can be chosen to be of the form ${\sf id}_{PA^*}$
for some finite prefix code $P$ with $\mu(P) = h$.

\noindent {\bf (2)} \ For all $g, h \in \mathbb{Z}[\frac{1}{k}]$
with \ $0 \leq g < h \leq 1$ \ the elements 
$\varphi_g, \varphi_h \in M_{k,1}$ constructed in {\bf (1)} satisfy
 \ \  $\varphi_g <_{\cal R} \varphi_h$.

Hence \ $\{\varphi_h : h \in \mathbb{Z}[\frac{1}{k}] \cap [0,1]\}$ \ forms
a dense $<_{\cal R}$-chain of elements of $M_{k,1}$ with 
${\sf height}_{\cal R}(\varphi_h) = h$.
\end{pro}
{\bf Proof.} {\bf (1)} We use Lemma \ref{prefcodesformeansure} (1) to 
construct a finite prefix code $P$ with \ $h = \mu(P)$. 
Letting $\varphi = {\sf id}_{PA^*}$ we obtain 
${\sf height}_{\cal R}({\sf id}_{PA^*}) = \mu(P)$.

\smallskip

\noindent {\bf (2)} We use Lemma \ref{prefcodesformeansure} (2), and let
$\varphi_h = {\sf id}_{P_hA^*}$.
Then, ${\sf id}_{P_hA^*} = {\sf id}_{P_hA^*} \ {\sf id}_{P_gA^*} = $
${\sf id}_{P_gA^*} \ {\sf id}_{P_hA^*}$, i.e., 
${\sf id}_{P_hA^*} \leq {\sf id}_{P_gA^*}$ \ (in the idempotent order). 
Moreover, since $\mu(P_h) < \mu(P_g)$ we have
 \ ${\sf id}_{P_hA^*} <_{\cal R} {\sf id}_{P_gA^*}$ \ and 
 \ ${\sf id}_{P_hA^*} <_{\cal L} {\sf id}_{P_gA^*}$.   
 \ \ \ $\Box$

%%%%%%%%%%%%%%%%%%%%%%%%%%%%%%

\subsection{A height function for the $\cal L$-order }

The construction of an $\cal L$-height (Definitions \ref{collision} and
\ref{Lheight} below) requires a few preliminary definitions and facts. 

Let $\simeq$ be a partition of a set $S \subseteq A^*$.
For each $x \in S$ the $\simeq$-class containing $x$ is denoted by $[x]$.
A set of {\it representatives} of $\simeq$ is, by definition, a set 
consisting of exactly one element from each $\simeq$-class.
A set of {\it minimum-length representatives} of $\simeq$ is a set $R$ of 
representatives such that each $r \in R$ has minimum length in $[r]$.
A set of maximum-length representatives is defined in a similar way.

\begin{lem} \label{classw_restrVSlength} \
For a right ideal homomorphism $\varphi$ let $P = {\sf domC}(\varphi)$
and assume that ${\sf part}(\varphi)$ is a prefix code congruence. Then we
have:

\smallskip

\noindent
{\bf (1)} \ Let $n$ be at least as large as the length of the longest word
in ${\sf imC}(\varphi)$.
There exists an essential restriction $\varphi_0$ of $\varphi$
such that all the words in ${\sf imC}(\varphi_0)$ have the same length,
i.e., ${\sf imC}(\varphi_0) \subseteq A^n$.

\smallskip

\noindent
{\bf (2)} \ There is an essential restriction $\varphi_1$ of $\varphi$ such
that the {\em minimum}-length representatives of the partition
$\simeq_{{\sf domC}(\varphi_1)}$ all have the same length, and 
${\sf part}(\varphi_1)$ is a prefix code congruence.

Similarly, there is an essential restriction $\varphi_2$ of $\varphi$ such
that all the {\em maximum}-length representatives of 
$\simeq_{{\sf domC}(\varphi_2)}$ have the same length, and
${\sf part}(\varphi_2)$ is a prefix code congruence.
\end{lem}
{\bf Proof.} (1) \ Let $P_i$ be a class of $\varphi$ in
${\sf domC}(\varphi)$, and let $\varphi(P_i) = y_i \in {\sf imC}(\varphi)$;
hence, $\varphi^{-1}(y_i) = P_i$.  We consider the essential restriction
that replaces the set of table entries
 \ $P_i \times \{y_i\} = \{(x,y_i) : x \in P_i\}$ \ by the set
 \ $\bigcup_{j=1}^k P_ia_j \times \{y_ia_j\} \ = \ $
$\{(xa_1,y_ia_1) : x \in P_i\}$ $ \ \cup \ \dots \ \cup \ $
$\{(xa_k,y_ia_k) : x \in P_i\}$.
By such a replacement we can make the shortest element of
${\sf imC}(\varphi)$ longer. By repeating this, we can give the same length
to all elements of  ${\sf imC}(\varphi)$.

(2) \ Similarly, a class-wise replacement step can be used to make all the 
elements of $P_i$ longer; in particular, minimum- (or maximum-) length 
elements can be made longer. 
By repeating this on the class that has the shortest among the minimum- 
(or maximum-) length elements over all classes, we can give the same
length to all the minimum-length (or maximum-length) representatives.
Since only class-wise replacements are used, ${\sf part}(\varphi_1)$ remains
a prefix code congruence.
   \ \ \ $\Box$

\medskip

\noindent {\bf Remark.} Any right ideal homomorphism $\varphi$ can be
essentially equally restricted to a right ideal homomorphism $\Phi$ such
that all words in ${\sf domC}(\Phi)$ have the same length. However, in 
general ${\sf part}(\Phi)$ will no longer necessarily be a prefix code 
congruence.

\medskip

For any finite prefix code $P \subset A^*$, a {\it complementary prefix code}
of $P$ is a finite prefix code $Q \subset A^*$ such that
${\sf ends}(PA^*) \cap {\sf ends}(QA^*) = \varnothing$
and ${\sf ends}(PA^*) \cup {\sf ends}(QA^*) = A^{\omega}$. 
This was introduced in Definition 3.29 in \cite{BiRL}. By Lemma 3.30 in 
\cite{BiRL}, every finite prefix code $P$ has a complementary prefix code
(which is empty iff $P$ is a maximal prefix code).

\medskip

We now start the construction of an $\cal L$-height function; this is more 
subtle than the $\cal R$-height since now we have to measure how fine
a partition is rather than just how large a set is.
Intuitively, elements $\varphi \in M_{k,1}$ that are higher in the
$\cal L$-order have ``smaller'' and ``more'' classes in 
${\sf part}(\varphi)$; here we should treat the complement of
${\sf domC}(\varphi)$ like a (virtual) class too (called the ``undefined
class''). 
All classes of ${\sf part}(\varphi)$ are finite (of size 
$\leq |{\sf domC}(\varphi)|$).
The highest elements in the $\cal L$-order (i.e., the $\cal L$-class of 
the identity and the injective total maps of $M_{k,1}$) only have 
singleton classes.
This suggests that the singleton classes should be given a largeness
of zero, and that the concept of ``collisions'' of a function is relevant 
for measuring the largeness of the classes.
A total injective function has no collisions.
In a non-injective function $f$, a {\it collision} is any pair $(x_1,x_2)$ 
such that $x_1 \neq x_2$ and $f(x_1) = f(x_2)$. The concept of collision 
is commonly used in Algorithms and Data Structures.
Thus, the first idea is to say that a class $f^{-1}(y)$ of $f$ has
$|f^{-1}(y)| - 1$ collisions, if $y \in {\sf Im}(f)$; the subtraction of 
$1$ is justified by the fact that one element by itself is not a collision
(collisions only start with the second element in a class).
Also, any $x$ for which $f(x)$ is undefined will be treated as a collision
all by itself; thus, the undefined class $C_{\varnothing}$ has
$|C_{\varnothing}|$ collisions.
Moreover, for $M_{k,1}$ we need to use the measure $\mu$ rather than
cardinality.

This motivates the following:

\begin{defn} \label{collision} \
Let $\simeq$ be a {\em prefix code congruence} on a right ideal $PA^*$, 
where $P$ is a finite prefix code. 
Let $\simeq_P$ be the restriction of $\simeq$ to $P$.
Let $\{P_1, \ldots, P_n\}$ be the classes of $\simeq_P$, and let
$P_{\varnothing}$ be a complementary prefix code of $P$ in $A^*$. 

The  amount of collision of $\simeq$ in $P_{\varnothing}$ is
$\mu(P_{\varnothing})$.
For a class $P_i$ of $\simeq_P$ ($1 \leq i \leq n$), let $m_i$ be any chosen 
minimum-length element in $P_i$. The amount of collision of $\simeq$ in 
$P_i$ is \ $\mu(P_i - \{ m_i \} )$.
The total {\bf amount of collision} of the prefix code congruence $\simeq$ 
is

\medskip

 \ \ \ \ \ ${\sf coll}(\simeq) \ = \ \mu(P_{\varnothing}) \ + \ $
$\sum_{i=1}^n \mu(P_i - \{m_i\}).$

\medskip

\noindent Since $P_{\varnothing} \cup P_1 \cup \ \ldots \ \cup P_n$ is a 
maximal prefix code, \ $\mu(P_{\varnothing}) + \sum_{i=1}^n \mu(P_i) = 1$; 
hence we also have

\medskip

 \ \ \ \ \ ${\sf coll}(\simeq) \ = \ 1 \ - \ \sum_{i=1}^n \mu(m_i). $

\medskip

\noindent Accordingly, the {\bf amount of non-collision} of the prefix 
code congruence $\simeq$ is defined by 

\medskip
 
 \ \ \ \ \ ${\sf noncoll}(\simeq) \ = \ \sum_{i=1}^n \mu(m_i). $
\end{defn}
Further justifications of this definition:

The motivation for removing an element $m_i$ from the class $P_i$ when we 
measure the collisions is that one element by itself creates no collision;
only subsequent additions of elements to a class cause collisions. We
choose to remove the most probable (i.e., shortest) element from each class, 
and let all the other elements in the class account for the collisions. 
At the end of this subsection there is a discussion of other possible 
definitional choices.

The value of ${\sf coll}(\simeq)$ depends only on $\simeq$. Indeed, first,
it is easy to see that ${\sf coll}(\simeq)$ does not depend on the choice 
of a particular minimum-length word $m_i$ in $P_i$, since 
${\sf coll}(\simeq)$ depends only on the lengths of words.
Second, we easily show that ${\sf coll}(\simeq)$ does not depend on the 
choice of $P_{\varnothing}$, since all complementary prefix codes of $P$ in 
$A^*$ have the same ends (namely $A^{\omega} - {\sf ends}(PA^*)$), hence 
the same measure.  
The values of ${\sf coll}(\simeq)$ are $k$-ary rational numbers 
that range from $0$ (for the identity congruence) to $1$ (for the empty 
congruence, on an empty domain).

\begin{lem} \label{coll_vs_essequ} \
If $\simeq'$ and $\simeq$ are prefix code congruences and
 \ $\simeq' \ =_{\sf ess} \ \simeq$ \ then
 \ ${\sf coll}(\simeq') = {\sf coll}(\simeq)$.
\end{lem}
{\bf Proof.} If we apply a class-wise replacement step
$C \to \{ Ca_1, \ldots, Ca_k\}$ to $\simeq$, 
where $C$ is a class of $\simeq$, the resulting prefix code
congruence $\simeq'$ satisfies:

\smallskip

${\sf coll}(\simeq') \ = \ {\sf coll}(\simeq) \ - \ \mu(C - \{m\}) $
 $ \ + \ \mu(Ca_1 - \{ma_1\}) + \ \ \ldots \ \ + \mu(Ca_k - \{ma_k\})$ , 

\smallskip

\noindent where $m$ is a minimum-length element of $C$. For any set $S$ 
we have \ $\mu(S) = \mu(Sa_1) + \ldots + \mu(Sa_k)$, since 
$\mu(Sa_i) = \frac{1}{k} \, \mu(S)$; it follows that
${\sf coll}(\simeq') = {\sf coll}(\simeq)$.

In a similar way one proves that an inverse class-wise replacement step 
preserves the amount of collision. Hence, iteration of replacement steps and
inverse replacement steps preserves the amount of collision.
 \ \ \ $\Box$

\medskip

More generally we have the following (note the order reversal, since finer
congruences have fewer collisions):

\begin{lem} \label{collOrder_vs_endOrder} \
Suppose $\simeq_1$ and $\simeq_2$ are prefix code congruences that are
comparable in the order $\leq_{\sf ends}$. Then we have
 \ $\simeq_1 \ <_{\sf ends} \ \simeq_2$, \ or
 \ $\simeq_1 \ =_{\sf ess}  \ \simeq_2$, \ or
 \ $\simeq_1 \ >_{\sf ends} \ \simeq_2$, \ according as
 \ ${\sf coll}(\simeq_1) > {\sf coll}(\simeq_2)$,\ or
 \ ${\sf coll}(\simeq_1) = {\sf coll}(\simeq_2)$,\ or
 \ ${\sf coll}(\simeq_1) < {\sf coll}(\simeq_2)$.
\end{lem}
{\bf Proof.} Suppose $\simeq_1 >_{\sf ends} \simeq_2$. Then, by Lemma
\ref{coll_vs_essequ}, we can essentially equally restrict $\simeq_1$ and 
$\simeq_2$ so that in the resulting prefix code congruences (which we still 
call $\simeq_1$ and $\simeq_2$) we have: Every class of $\simeq_2$ is a 
union of classes of $\simeq_1$.

Suppose  $Q$ is a class of $\simeq_2$ in ${\sf domC}(\simeq_2)$, and suppose
$P_1, \ldots, P_n$ are the classes of $\simeq_1$ in ${\sf domC}(\simeq_1)$
such that $Q = P_1 \cup \ldots \cup P_n$, \ $2 \leq n$.
Then the amount of collision in $\simeq_2$ for $Q$ is $\mu(Q - m)$ (where
$m$ is a shortest element of $Q$).
The amount of collision in $\simeq_1$ for $P_1, \ldots, P_n$ (with shortest
element in $P_i$ denoted by $m_i$) is

\smallskip

$\mu(P_1 - m_1) + \ \ldots \ + \mu(P_n - m_n) \ = \ $
$\mu(P_1) + \ \ldots \ + \mu(P_n) \ - \ \mu(m_1) - \ \ldots \ - \mu(m_n)$

\smallskip

$ \ = \ $ $\mu(Q) \ - \ \mu(m_1) - \ \ldots \ - \mu(m_n) \ < \ $
$\mu(Q) \ - \ \mu(m)$.

\smallskip

\noindent The last ``$<$'' is due to the fact that $\mu(m)$ is equal to one
of the numbers $\mu(m_1), \ldots, \mu(m_n)$, since
$Q = P_1 \cup \ \ldots \ \cup P_n$ and $n \geq 2$. We conclude that

\smallskip

$\mu(P_1 - m_1) + \ \ldots \ + \mu(P_n - m_n) \ < \ \mu(Q - m)$.

\smallskip

\noindent In other words, coarser classes have larger amounts of collision.

Moreover, if $C$ is a class of $\simeq_1$ that does not intersect the 
domain of $\simeq_2$, then $C$ is in the undefined class of $\simeq_2$, 
hence $\mu(C)$ will be counted in the amount collision in $\simeq_2$ (but 
only $\mu(C -m)$ will be counted in $\simeq_1$). So, 
here again, the amount of collision in $\simeq_2$ is larger.
So, ${\sf coll}(\simeq_1) < {\sf coll}(\simeq_2)$.

\smallskip

In a  similar way we can prove that $\simeq_1 \ <_{\sf ends} \ \simeq_2$ 
implies ${\sf coll}(\simeq_1) > {\sf coll}(\simeq_2)$. And the proof
that $\simeq_1 \ =_{\sf ess}  \ \simeq_2$ implies
${\sf coll}(\simeq_1) = {\sf coll}(\simeq_2)$ was already given in
Lemma \ref{coll_vs_essequ}.

\smallskip

For the converse: 
Suppose we have ${\sf coll}(\simeq_1) > {\sf coll}(\simeq_2)$ and suppose
that $\simeq_1$ and $\simeq_2$ are comparable for the 
$\leq_{\sf ends}$-order. 
This leaves only the three possibilities: \ $<_{\sf ends}$, $=_{\sf ess}$, 
and $>_{\sf ends}$. But we already proved that $=_{\sf ess}$ and 
$>_{\sf ends}$ would contradict 
${\sf coll}(\simeq_1) > {\sf coll}(\simeq_2)$. So we
have $\simeq_1 \ <_{\sf ends} \ \simeq_2$.
 \ \ \ \ \ $\Box$

\begin{defn} \label{Lheight} \
For any element of $M_{k,1}$ represented by a right ideal homomorphism 
$\varphi$ we define the $\cal L$-height by
 \ ${\sf height}_{\cal L}(\varphi) \ = \ 1 - {\sf coll}({\sf part}(\varphi))$
 \ (i.e., the amount of non-collision). Hence, 

\medskip

\hspace{1.0in}
${\sf height}_{\cal L}(\varphi) \ = \ \sum_{i=1}^n \mu(m_i)$ ,

\medskip

\noindent
where $m_i$ is a shortest representative of the class $P_i$ of
${\sf part}(\varphi)$ (denoting the classes of ${\sf part}(\varphi)$ in 
${\sf domC}(\varphi)$ by $P_1, \ldots, P_n$).
\end{defn}
By the characterization of the $\cal L$-order of $M_{k,1}$ and by Lemma 
\ref{coll_vs_essequ} above,  ${\sf height}_{\cal L}(\varphi)$ depends only 
on $\varphi$ as an element of $M_{k,1}$ and not on the right ideal 
homomorphism chosen.
Lemma \ref{collOrder_vs_endOrder} implies that ${\sf height}_{\cal L}(.)$ 
is indeed a height function for $\leq_{\cal L}$, i.e., that we have: 

\begin{pro} \label{Lheight_vs_Lorder} \
Suppose $\varphi, \psi \in M_{k,1}$ are comparable in the $\cal L$-order.
Then we have \ $\varphi >_{\cal L} \psi$, \ or
 \ $\varphi \equiv_{\cal L} \psi$, \ or
 \ $\varphi <_{\cal L} \psi$, \ according as
 \ ${\sf height}_{\cal L}(\varphi) > {\sf height}_{\cal L}(\psi)$,  \ or
 \ ${\sf height}_{\cal L}(\varphi) = {\sf height}_{\cal L}(\psi)$,  \ or
 \ ${\sf height}_{\cal L}(\varphi) < {\sf height}_{\cal L}(\psi)$.
 \ \ \ \ \ $\Box$
\end{pro}

\begin{pro} \label{collOrder_forh} \
{\bf (1)} \ For every $h \in \mathbb{Z}[\frac{1}{k}] \cap [0,1]$ 
there exists $\varphi_h \in M_{k,1}$ such that 
 \ ${\sf height}_{\cal L}(\varphi_h) = h$.

\noindent {\bf (2)} \ For all $g, h \in \mathbb{Z}[\frac{1}{k}]$
with \ $0 \leq g < h \leq 1$ \ the elements
$\varphi_g, \varphi_h \in M_{k,1}$ constructed in {\bf (1)} satisfy
 \ \  $\varphi_g <_{\cal L} \varphi_h$.
The set \ $\{\varphi_h : h \in \mathbb{Z}[\frac{1}{k}] \cap [0,1]\}$ \ forms
a dense $<_{\cal L}$-chain of elements of $M_{k,1}$ with 
${\sf height}_{\cal L}(\varphi_h) = h$.
\end{pro}
{\bf Proof.} This is proved in the same way as 
Prop.\ \ref{every-rat-is_Rheight}. 
 \ \ \ $\Box$

%\bigskip

%\bigskip

\newpage 

\noindent 
{\bf Variants of the definition of an $\cal L$-height function:}

\medskip

We chose the definition 
 \ ${\sf height}_{\cal L}(\varphi) = \sum_{i=1}^n \mu(m_i)$
 \ where each word $m_i$ is a {\it minimum-length} representative of a
${\sf part}(\varphi)$-class in ${\sf domC}(\varphi)$. For the remainder of
this subsection we will call this function 
${\sf height}_{\cal L}^{\sf min}(.)$.
If in the definition of $\cal L$-height we replace minimum-length by
{\it maximum-length} representatives we obtain a function 
${\sf height}_{\cal L}^{\sf max}(.)$ that is also an $\cal L$-height 
function (according to Def.\ \ref{Heightfunct}).
For all $\varphi \in M_{k,1}$ we obviously have 
 \ ${\sf height}_{\cal L}^{\sf max}(\varphi) \leq $
${\sf height}_{\cal L}^{\sf min}(\varphi)$.
For idempotents the following relation holds between the $\cal R$-height
function and the two $\cal L$-height functions.

\begin{pro} \label{RvsLheightsIdempotents} \
For any idempotent $\eta = \eta^2 \in M_{k,1}$,  

\medskip

\hspace{1.in}
${\sf height}_{\cal L}^{\sf max}(\eta) \ \leq \ {\sf height}_{\cal R}(\eta)$
$ \ \leq \ {\sf height}_{\cal L}^{\sf min}(\eta)$.

\medskip

\noindent If $\eta = \eta^2 \in {\it Inv}_{k,1}$ then 
 \ ${\sf height}_{\cal L}^{\sf max}(\eta) \ = \ $
${\sf height}_{\cal R}(\eta) \ = \ $
${\sf height}_{\cal L}^{\sf min}(\eta)$.
\end{pro}
{\bf Proof.} For an idempotent, the elements of ${\sf imC}(\eta)$ form a
set of representatives of the ${\sf part}(\eta)$-classes in
${\sf domC}(\eta)$, assuming that $\eta$ has been restricted so as to make
${\sf part}(\eta)$ a prefix code congruence. Hence, the lengths of the elements
of ${\sf imC}(\eta)$ are between the lengths of the minimum-length
representatives and the maximum-length representatives. The inequalities
follow.  

When the idempotent $\eta$ is injective, the congruence classes of 
${\sf part}(\eta)$ are singletons, so the minimum-length and the 
maximum-length representatives of a class are the same.  
 \ \ \ $\Box$

\medskip

\noindent An $\cal L$-height function could also be defined by using the
{\em average} of the lengths in each block of ${\sf part}(\varphi)$:

\smallskip

\hspace{1.in} ${\sf height}_{\cal L}^{\sf ave}(\varphi) \ = \ $
$\sum_{y \in {\sf imC}(\varphi)} \ \mu( {\sf ave}(\varphi^{-1}(y)) )$ , 

\smallskip

\noindent where for any finite set $S \subset A^*$ we define 
 \ ${\sf ave}(S) = \frac{1}{|S|} \sum_{x \in S} |x|$ .  
This is indeed an  $\cal L$-height function, as a consequence of the 
fact that for two disjoint finite sets $S_1, S_2$ we have

\smallskip

\hspace{1.in}
 ${\sf min}\{ {\sf ave}(S_1), {\sf ave}(S_2)\} \ \leq \ $
${\sf ave}(S_1 \cup S_2) \ \leq \ $
$ {\sf max}\{ {\sf ave}(S_1), {\sf ave}(S_2)\}$.

\smallskip

\noindent
The same reasoning works with the average replaced by the {\em median}.

%%%%%%%%%%%%%%%%%%%%%%%%%%%%%%%

\subsection{A connection between the $\cal R$- and $\cal L$-heights and
the $\cal D$-relation }

Interestingly, ${\sf height}_{\cal R}(\varphi)$ determines the 
$\cal D$-class of $\varphi$; similarly, ${\sf height}_{\cal L}(\varphi)$ 
determines the $\cal D$-class. First of all, obviously, 

\smallskip

 \ \ \ \ \  ${\sf height}_{\cal L}(\varphi) = 0$ \ \ iff
 \ \ $\varphi = {\bf 0}$ \ \ iff 
 \ \ ${\sf height}_{\cal R}(\varphi) = 0$.

\smallskip

\noindent When $\varphi \neq {\bf 0}$ it was proved (Theorem 2.5 in 
\cite{BiThomMon}) that the $\cal D$-class of $\varphi$ is 

\smallskip

 \ \ \ \ \ \  $D(\varphi) \ = \  $ 
$\{ \psi \in M_{k,1} : \ |{\sf imC}(\psi)| \equiv |{\sf imC}(\varphi)|$ 
 \ {\sf mod} $k-1 \}$ .

\smallskip

\noindent The number $i \in \{1, \ldots, k-1\}$ such that 
 \ $i \equiv |{\sf imC}(\varphi)|$ {\sf mod} $k-1$ \ is called the 
{\it index of the $\cal D$-class of} $\varphi$ (when $\varphi \neq {\bf 0}$).
In this paper, integers modulo $k-1$ will be picked in the range
$\{1, \ldots, k-1\}$.

Recall that for a $k$-ary rational number $r = a/k^n$, where $k$ does not 
divide $a$, the numerator $a$ is denoted by ${\sf num}(r)$.

\begin{pro} \label{heightdetDclass} \    
For every $\varphi \in M_{k,1}$ the $\equiv _{\cal D}$-class of $\varphi$
is uniquely determined by ${\sf height}_{\cal R}(\varphi)$
(and similarly, by ${\sf height}_{\cal L}(\varphi)$). \ More precisely, 
when $\varphi \neq {\bf 0}$ we have the following formulas.

\medskip
 
\noindent {\bf (1)} \ The $\cal D$-class index $i \in \{1, \ldots, k-1\}$ 
of $\varphi$ is determined by

\smallskip

 \ \ \   \ \ \ \ \   
$i \ \equiv \ {\sf num}({\sf height}_{\cal L}(\varphi)) \ \equiv \ $
${\sf num}({\sf height}_{\cal R}(\varphi))$ \ \ {\sf mod} \ $k-1$.

\medskip

\noindent {\bf (2)} \ When base-$k$ representations of 
${\sf height}_{\cal L}(\varphi)$ and ${\sf height}_{\cal R}(\varphi)$ are 
given we have:

\smallskip

$\bullet$ \ \ If 
 \ ${\sf height}_{\cal L}(\varphi) \ = \ 0.d_1 \ldots d_m$ \ or 
 \ ${\sf height}_{\cal R}(\varphi) \ = \ 0.d'_1 \ldots d'_n$ \ we have
 
\smallskip

 \ \ \ \ \   \ \ \    
$i \ \equiv \ d_1 + \ldots + d_m \ \ {\sf mod} \ k-1$ , \ \ or  

\smallskip

 \ \ \ \ \  \ \ \   
$i \ \equiv \ d'_1 + \ldots + d'_n \ \ {\sf mod} \ k-1$ .

\smallskip

$\bullet$ \ \ If ${\sf height}_{\cal L}(\varphi) = 1$ or
${\sf height}_{\cal R}(\varphi) = 1$ then $i = 1$. 
\end{pro}
{\bf Proof.} \ Part (2) immediately follows from part (1), since
$k \equiv 1$ {\sf mod} $k-1$. Let us prove part (1). 

\smallskip

$\cal R$-height formula: \ We have 
${\sf height}_{\cal R}(\varphi) = \mu({\sf imC}(\varphi))$, and
we can write $|{\sf imC}(\varphi)| = i + j \, (k-1)$, for some integers
$i, j$ such that $1 \leq i \leq k-1$ and $j  \geq 0$. Moreover,
$\mu({\sf imC}(\varphi)) = (i + j \, (k-1)) \cdot k^{-N}$, for some 
integer $N > 0$. 
Note that for the {\sf mod} $k-1$ value of the numerator of a $k$-ary 
rational number, it does not matter whether the numerator is divisible
by $k$ (since $k \equiv 1$ {\sf mod} $k-1$). Hence,
$|{\sf imC}(\varphi)| \ \equiv \ i \ \equiv \ $
${\sf num}\big(\mu({\sf imC}(\varphi))\big)$ \ {\sf mod} $k-1$.

\smallskip

$\cal L$-height formula: \ We have  
${\sf height}_{\cal L}(\varphi) = \sum_{i=1}^n \mu(m_i)$, where 
$n = |{\sf imC}(\varphi)|$, and $\{m_i: i = 1, \ldots,n\}$ is the set of 
minimum-length representatives of the ${\sf part}(\varphi)$ classes in
${\sf domC}(\varphi)$.
By Lemma \ref{classw_restrVSlength} we can assume that all $m_i$ have the 
same length, say $\ell$. Then, 
${\sf height}_{\cal L}(\varphi) = n \, k^{-\ell}$. 
Again, for the {\sf mod} $k-1$ value of the numerator it does not matter
whether the numerator is divisible by $k$. Hence, 
$|{\sf imC}(\varphi)| \ = \ n \ \equiv \ $
${\sf num}({\sf height}_{\cal L}(\varphi))$ \ {\sf mod} $k-1$. 
 \ \ \ $\Box$

\begin{pro} \label{indepLhRh} 
{\bf (Independence of the $\cal R$- and $\cal L$-heights in $M_{k,1}$ and 
in ${\it Inv}_{k,1}$).} \\     
Let $h_1, h_2$ be any $k$-ary rationals with $0 < h_1, h_2 \leq 1$ and 
such that ${\sf num}(h_1) \equiv {\sf num}(h_2)$ {\sf mod} $k-1$ (i.e., 
$h_1$ and $h_2$ determine the same non-zero $\equiv_{\cal D}$-class).
Then there exists an element $\varphi \in {\it Inv}_{k,1}$ 
$(\subset M_{k,1})$ such that 
 \ ${\sf height}_{\cal L}(\varphi) = h_1$ \ and 
 \ ${\sf height}_{\cal R}(\varphi) = h_2$. 
\end{pro}
{\bf Proof.} This follows directly from Lemma \ref{LemmaindepLhRh}, which 
will be proved next. \ \ \ $\Box$

\begin{lem}  \label{LemmaindepLhRh} . \\  
{\bf (1)} \ Let $R$ be any non-zero $\cal R$-class of $M_{k,1}$ and let
$h_1$ be any $k$-ary rational with $0 < h_1 \leq 1$, such that the
$\cal D$-class of $R$ coincides with the $\cal D$-class determined by $h_1$
(used as an $\cal L$-height); in other words, we assume that
 \ ${\sf num}({\sf height}_{\cal R}(R)) \ \equiv \ {\sf num}(h_1)$
 \ {\sf mod} $k-1$.
Then there exists an element $\varphi \in R \cap {\it Inv}_{k,1}$
such that \ ${\sf height}_{\cal L}(\varphi) = h_1$.

\smallskip

\noindent {\bf (2)} \ Similarly, let $L$ be any non-zero $\cal L$-class and
let $h_2$ be any $k$-ary rational with $0 < h_2 \leq 1$, such that the
$\cal D$-class of $L$ coincides with the $\cal D$-class determined by $h_2$
(used as an $\cal R$-height).
Then there exists an element $\psi \in L \cap {\it Inv}_{k,1}$ such that
 \ ${\sf height}_{\cal R}(\psi) = h_2$.
\end{lem}
{\bf Proof.} We only prove part (1), since (2) is similar.
We consider the base-$k$ representation 
 \ $h_1 = 0.d_1 \ldots d_n$ \ (or $h_1 = 1$).
As in the proof of Lemma \ref{prefcodesformeansure}, we construct the 
finite prefix code $P_{h_1}$ from $h_1$, with $\mu(P_{h_1}) = h_1$
and \ $|P_{h_1}| = \sum_{i=1}^n d_i$ \ (or $\mu(P_{h_1}) = |P_{h_1}| = 1$ 
if $h_1 = 1$).
We can increase the size of $P_{h_1}$ by multiples of $k-1$, as follows:
 \ In $P_{h_1}$ we replace \ $a_{d_1+1} \ldots a_{d_{n-1}+1} a_{d_n}$ 
 \ by \ $a_{d_1+1} \ldots a_{d_{n-1}+1} a_{d_n} \, A$, thus obtaining
a prefix code $P_{h_1}^{(1)}$ of size $|P_{h_1}| + k-1$, with the measure 
remaining unchanged at $\mu(P_{h_1}^{(1)}) = h_1$ (by Lemma
\ref{set_rewrite}).
This can be repeated: In $P_{h_1}^{(1)}$ we replace 
 \ $a_{d_1+1} \ldots a_{d_{n-1}+1} a_{d_n} a_k$ \ by
 \ $a_{d_1+1} \ldots a_{d_{n-1}+1} a_{d_n} a_k \, A$, thus obtaining
a prefix code $P_{h_1}^{(2)}$ of size \ $|P_{h_1}| + 2 \cdot (k-1)$, with
unchanged measure $\mu(P_{h_1}^{(2)}) = h_1$.
As a result, for any $j \geq 1$ we obtain a prefix code 
$P_{h_1}^{(j)}$ of size \ $|P_{h_1}| + j \cdot (k-1)$, with
unchanged measure $\mu(P_{h_1}^{(j)}) = h_1$.
More precisely, inductively, $P_{h_1}^{(j)} = $
$(P_{h_1}^{(j-1)} - \{ a_{d_1} \ldots a_{d_{n-1}+1} a_k^j\})$  $ \cup $
$ a_{d_1} \ldots a_{d_{n-1}+1}  a_k^j A$

For all $\psi$ in the $\cal R$-class $R$ the right ideals ${\sf Im}(\psi)$ 
are essentially equal, and all are essentially equal to $QA^*$ for a fixed 
finite prefix 
code $Q$. Since $R$ and $h_1$ correspond to the same $\cal D$-class, we
have \ $|Q| \equiv {\sf num}(h_1)$ {\sf mod} $k-1$. 
We can increase the size of $Q$ by any multiple of $k-1$ without
changing the corresponding $\cal R$-class, as follows: For any $q \in Q$
we replace $Q$ by \ $Q^{(1)} = (Q - \{q\}) \, \cup \, q A$, of size 
$|Q^{(1)}| = |Q| + k-1$, and such that \ $Q^{(1)}A^* =_{\sf ess} QA^*$. 
Then we replace $Q^{(1)}$ by 
 \ $Q^{(2)}  = (Q^{(1)} - \{qa_k\}) \, \cup \, qa_kA$, etc. After $j$ 
steps we obtain a prefix code $Q^{(j)}$ of size 
 \ $|Q^{(j)}| = |Q| + j \cdot (k-1)$, 
such that \ $Q^{(j)} A^* =_{\sf ess} QA^*$.

Since $R$ and $h_1$ correspond to the same $\cal D$-class, i.e.,
$|Q| \equiv {\sf num}(h_1) = |P_{h_1}|$ {\sf mod} $k-1$, there exist
$j$ and $j'$ such that \ $|P_{h_1}^{(j)}| = |Q^{(j')}|$.
Let us define $\varphi$ by any bijection \ $P_{h_1}^{(j)} \to Q^{(j')}$. 
Then $\varphi \in R$ since \ $Q^{(j')}A^* =_{\sf ess} QA^*$.
Also, ${\sf height}_{\cal L}(\varphi) = h_1$ \ since $\varphi$ is 
injective and since ${\sf domC}(\varphi) = P_{h_1}^{(j)}$ with 
$\mu(P_{h_1}^{(j)}) = h_1$.
 \ \ \ $\Box$

%\subsection{ $\cal R$- versus $\cal L$-heights of idempotents }

\bigskip

We saw in Prop.\ \ref{RvsLheightsIdempotents} that for idempotents of 
$M_{k,1}$, there are relations between the $\cal R$-height and the
$\cal L$-height.

%%%%%%%%%%%%%%%%%%%%%%%%%%%%%%%%%%%%%%%%%%%%%%%%%%%
%% Section 
%%%%%%%%%%%%%%%%%%%%%%%%%%%%%%%%%%%%%%%%%%%%%%%%%%%

%%%%%

\section{The Green relations of ${\it plep}M_{k,1}$ and 
         ${\it tlep}M_{k,1}$ }

%%%

\subsection{ The monoids ${\it plep}M_{k,1}$ and ${\it tlep}M_{k,1}$ }

\noindent
The submonoid ${\it tlep}M_{k,1}$ of {\it total length-equality preserving} 
elements of $M_{k,1}$ was introduced in \cite{BiThomMon}, where it was
simply called ${\it lep}M_{k,1}$. 
We now add the ``{\it t}'' (for total) in order to distinguish it from the 
submonoid ${\it plep}M_{k,1}$ of {\it partial length-equality preserving} 
elements of $M_{k,1}$. As usual, partial does not rule out total, so 
 \ ${\it tlep}M_{k,1} \subset {\it plep}M_{k,1}$. 
More precisely, these submonoids of $M_{k,1}$ are defined as follows.

\smallskip

${\it plep}M_{k,1} \ = \ \ \{ \varphi \in M_{k,1} : \ $
for all $x_1, x_2 \in {\sf Dom}(\varphi), \ |x_1| = |x_2|$
implies $|\varphi(x_1)| = |\varphi(x_2)| \}$.

\medskip

${\it tlep}M_{k,1} \ = \ \ \{ \varphi \in {\it plep}M_{k,1} : \ $
  ${\sf Dom}(\varphi)$ is an essential right ideal\}.

\medskip

\noindent In words, $\varphi \in M_{k,1}$ belongs to ${\it plep}M_{k,1}$ 
iff $\varphi$ transforms equal-length inputs to equal-length outputs, 
hence the name ``length equality preserving''. Recall that 
${\sf Dom}(\varphi)$ is essential iff 
${\sf ends}({\sf Dom}(\varphi)) = A^{\omega}$, i.e., iff
$\varphi$ is total on $A^{\omega}$.  

One can easily prove the following characterization: 
 \ $\varphi \in M_{k,1}$ belongs to ${\it plep}M_{k,1}$ \ iff 
 \ there is an essentially equal restriction $\Phi$ of $\varphi$ such that 
for some $m,n > 0$,

\smallskip

 \ \ \ ${\sf domC}(\Phi) \subseteq A^m$ \ and 
  \ ${\sf imC}(\Phi) \subseteq A^n$.  

\smallskip

\noindent For ${\it tlep}M_{k,1}$ we have in
addition that \ ${\sf domC}(\Phi) = A^m$.

\medskip

An important motivation for the study of ${\it plep}M_{k,1}$ and
${\it tlep}M_{k,1}$ is their similarity to (partial) acyclic boolean 
circuits.  In \cite{BiDistor} it was proved that ${\it tlep}M_{k,1}$ has a 
generating set of the form $\Gamma \cup \tau$ where $\Gamma$ is finite and 
$\tau = \{ \tau_{i,i+1} : i \geq 1\}$. Each $\tau_{i,i+1}$ is a position
transposition (or ``wire crossing''), defined as follows:
 \ $\tau_{i,i+1}(uabv) = ubav$ \ 
for all $u \in A^{i-1}$, \ $a, b \in A$, \ $v \in A^*$; and 
$\tau_{i,i+1}(x)$ is undefined when $|x| < i+1$.
When $k=2$, the set $\Gamma$ can be chosen to be 
 \ $\{{\sf and}, {\sf or}, {\sf not}, {\sf fork} \}$.
These are the classical circuit gates, given by the tables
 \ ${\sf and} = \{(00,0), (01,0), (10,0), (11,1)\}$, 
 \ ${\sf or}  = \{(00,0), (01,1), (10,1), (11,1)\}$, 
 \ ${\sf not} = \{(0,1), (1,0)\}$, 
 \ ${\sf fork} = \{(0,00), (1,11)\}$. 

It was proved in \cite{BiDistor} that for elements in ${\it tlep}M_{2,1}$, 
word-length over $\Gamma \cup \tau$ is polynomially equivalent to circuit-size.
For this reason we call generating sets of ${\it tlep}M_{k,1}$ (or, more
generally, of $M_{k,1}$, or of ${\it plep}M_{k,1}$) of the form 
$\Gamma \cup \tau$ (where $\Gamma$ is finite and $\tau$ is as above) 
{\it circuit-like generating sets}. 
The monoids $M_{k,1}$ and ${\it tlep}M_{k,1}$ have circuit-like generating 
sets, and Prop.\ \ref{cirlikeGenplepM} below will show the same for  
${\it plep}M_{k,1}$.  

If a different $\Gamma$ is used, the word-length changes only linearly,
by the following general observation (whose proof is straightforward):

\begin{pro} \label{finChangeGenSet} \ 
If two (possibly infinite) generating sets $\Gamma_1$ and $\Gamma_2$ for 
a monoid $M$ differ only by a finite amount (i.e., their symmetric 
difference $\Gamma_1 \! \vartriangle \! \Gamma_2$ is finite),
then the word-lengths of $M$ over $\Gamma_1$, respectively $\Gamma_2$,
are linearly related.
\ \ \ \ \  $\Box$
\end{pro} 
\noindent {\sf Notation:}  
When $P$ is a prefix code we abbreviate the partial identity map
${\sf id}_{PA^*}$ by ${\sf id}_P$. E.g., denoting the elements of the 
alphabet $A$ by $\{a_1, \ldots, a_k\}$, the partial identity 
${\sf id}_{A - a_1}$ is undefined on words that start with $a_1$ and is 
the identity on all other words in $A^*$.

\begin{pro} \label{cirlikeGenplepM} \  
The monoid ${\it plep}M_{k,1}$ has a circuit-like generating set.
More specifically, if $\Gamma \cup \tau$ is any circuit-like generating set
of ${\it tlep}M_{k,1}$ then 
 \ $\Gamma \cup \tau \cup \{ {\sf id}_{A - a_1} \}$ 
 \ generates ${\it plep}M_{k,1}$.
\end{pro} 
{\bf Proof.} \ For any $\varphi \in {\it plep}M_{k,1}$ with 
${\sf domC}(\varphi) \subseteq A^m$, we can define an element 
$\psi \in {\it tlep}M_{k,1}$ by extending the domain of $\varphi$ as 
follows: \ if $x \in {\sf domC}(\varphi)$ then $\psi(x) = \varphi(x)$,
and if $x \in A^m - {\sf domC}(\varphi)$ then $\psi(x) = y_0$,
where $y_0$ is any fixed element chosen in ${\sf imC}(\varphi)$.
Then we have: \ $\varphi = \psi \circ {\sf id}_{{\sf domC}(\varphi)}$.
This shows that ${\it plep}M_{k,1}$ is generated by $\Gamma \cup \tau$,
together with the partial identities of the form ${\sf id}_P$ (where $P$
ranges over the finite subsets of $A^m$ for all non-negative integers $m$).
Moreover, for $P \subseteq A^m$ we have 
$$ {\sf id}_P \ = \prod_{s \in A^m - P} {\sf id}_{A^m - \{s\}} $$
(and this composition of partial identities is commutative).
So, it will suffice to prove that each partial identity of the form
 \ ${\sf id}_{A^m - \{s\}}$ \ is generated by 
 \ $\Gamma \cup \{ {\sf id}_{A - a_1} \} \cup \tau$ \ for some finite
subset $\Gamma$ of ${\it tlep}M_{k,1}$.

For each letter $a_i \in A = \{a_1, a_2, \ldots, a_k\}$ we introduce the 
function $E_{a_i}: A \to \{a_1,a_2\}$, defined by \ $E_{a_i}(a_j) = a_1$ 
if $a_j \neq a_i$, and $E_{a_i}(a_j) = a_2$ if $a_j = a_i$. 
We also introduce the function \ ${\sf and}: A^2 \to \{a_1,a_2\}$, defined by 
 \ ${\sf and}(a_i a_j) = a_1$ \ if $a_i = a_1$ or $a_j = a_1$, and 
 \ ${\sf and}(a_i a_j) = a_2$ \ if $a_i \neq a_1 \neq a_j$.
And we define \ ${\sf not}: A \to \{a_1,a_2\}$ \ by ${\sf not}(a_1) = a_2$,
and ${\sf not}(a_i) = a_1$ for $a_i \in A - \{a_1\}$. 
Thus, the letter $a_1$ plays the role of the boolean value {\sf false}, and
the other letters play the role of {\sf true}.
We also use the function \ ${\sf fork}: A \to A^2$, defined by
 \ ${\sf fork}(a) = a \, a$. And we use \ ${\sf proj}_2: A^2 \to A$, defined 
by \ ${\sf proj}_2(a_i a_j) = a_j$.
Then \ ${\sf id}_{A^m - \{s\}}$ \ (for any $s = s_1 \ldots s_m \in A^m$) is 
generated by 

\smallskip

$\{ E_{a_i} : i = 1, \ldots, k\} \ \ \cup \ \ $
$\{ {\sf and}, \ {\sf not}, \ {\sf fork}, \ {\sf proj}_2,$ 
$ \ {\sf id}_{A - a_1} \} \ \ \cup \ \ \tau$.  

\smallskip 

\noindent Let us show how to simulate ${\sf id}_{A^m - \{s\}}$ 
by a fixed sequence of elements of the above set.  On input 
$x_1 x_2 \ldots x_m \ \ x_1 x_2 \ldots x_m \ w$ (where $x_i \in A$ for 
$i = 1, \ldots, m$ and $w \in A^*$), the output should be either undefined 
(if $x_1 x_2 \ldots x_m = s_1 s_2 \ldots s_m$),
or equal to the input (if $x_1 x_2 \ldots x_m \neq s_1 s_2 \ldots s_m$).

\smallskip

\noindent
First, by using $m$ copies of the fork function, together with 
transpositions $(\in \tau)$, a second copy of $x_1 x_2 \ldots x_m$ is
made:

\smallskip 

$x_1 x_2 \ldots x_m \ w \ \ \longmapsto \ \ $
 $x_1 x_2 \ldots x_m \ \ x_1 x_2 \ldots x_m \ w$ . 

\smallskip 

\noindent 
Next, by using $E_{s_1}, E_{s_2}, \ldots, E_{s_m}$ and transpositions,
we implement 

\smallskip 
 
$x_1 x_2 \ldots x_m \ \ x_1 x_2 \ldots x_m \ w \ \ \longmapsto \ \ $
$e_1 e_2 \ldots e_m \ \ x_1 x_2 \ldots x_m \ w$,

\smallskip 

\noindent where $e_i = a_1$ if $x_i \neq s_i$, and $e_i = a_2$ if 
$x_i = s_i$ ($1 \leq i \leq m$).
Then, to \ $e_1 e_2 \ldots e_m$ \ we apply $m-1$ copies of {\sf and}, as 
well as transpositions; this is followed by one application of {\sf not}; 
in effect, we compute the $m$-input {\sf nand} of $e_1 e_2 \ldots e_m$. 
This implements

\smallskip 

$e_1 e_2 \ldots e_m \ \ x_1 x_2 \ldots x_m \ w \ \ \longmapsto \ \ $
$e \ x_1 x_2 \ldots x_m \ w$,

\smallskip 

\noindent where $e = a_1$ if $e_1 = e_2 = \ \ldots \ = e_m = a_1$ \ (i.e.,
if \ $x_1 x_2 \ldots x_m = s_1 s_2 \ldots s_m$), and $e = a_2$ otherwise.
We apply ${\sf id}_{A - a_1}$ now; the operation is undefined if $e = a_1$, 
and is the identity otherwise. 
Finally, applying ${\sf proj}_2$ produces the output
 \ $x_1 x_2 \ldots x_m \ w$ \ if $e$ was $a_2$ (i.e., if 
${\sf id}_{A - a_1}$ was defined); the result is undefined otherwise.
 \ \ \ $\Box$

\bigskip

\noindent {\bf Open problems:} Are ${\it plep}M_{k,1}$ and
${\it tlep}M_{k,1}$ (not) finitely generated?
Is $M_{k,1}$ (not) finitely presented?

%%%

\subsection{ The $\cal R$-, $\cal L$-, and $\cal J$-relations of
${\it tlep}M_{k,1}$ and ${\it plep}M_{k,1}$ }

We will show that the $\cal R$-, $\cal L$-, and $\cal J$-orders of
${\it plep}M_{k,1}$ and ${\it tlep}M_{k,1}$ are very similar to the ones 
in $M_{k,1}$, and that ${\it plep}M_{k,1}$ is also congruence-simple.  

The monoids ${\it plep}M_{k,1}$ and ${\it tlep}M_{k,1}$ are {\it regular};
this is easily proved from the definition.
(A semigroup $S$ is called regular iff for every $s \in S$ there exists
$t \in S$ such that $sts = s$.)

\begin{pro} \label{L_R_plepM} \
The $\cal R$- and $\cal L$-orders of ${\it plep}M_{k,1}$ and of 
${\it tlep}M_{k,1}$  are induced by the corresponding orders of 
$M_{k,1}$. \ In other words, for all  
$\varphi_1, \varphi_2 \in {\it plep}M_{k,1}$ , 

\smallskip
 
 \ $\varphi_1 \geq_{\cal R} \varphi_2$ in ${\it plep}M_{k,1}$ \ iff
 \ $\varphi_1 \geq_{\cal R} \varphi_2$ in $M_{k,1}$ ; 

\smallskip

 \ $\varphi_1 \geq_{\cal L} \varphi_2$ in ${\it plep}M_{k,1}$ \ iff
 \ $\varphi_1 \geq_{\cal L} \varphi_2$ in $M_{k,1}$ . 

\smallskip

\noindent The $\cal R$- and $\cal L$-orders of ${\it tlep}M_{k,1}$ 
are also induced by the corresponding orders of $M_{k,1}$.
\end{pro}
{\bf Proof.} \ This follows from the fact that ${\it plep}M_{k,1}$ and
${\it tlep}M_{k,1}$ are regular semigroups, and the general fact that if 
$S_2$ is a subsemigroup of a semigroup $S_1$ and $S_2$ is regular then the 
$\geq_{\cal R}$ and $\geq_{\cal L}$ orders of $S_2$ are induced from 
$S_1$. See e.g.\ \cite{CliffPres} or p.\ 289 of \cite{Grillet}.  
 \ \ \ $\Box$

\medskip

An immediate consequence of Prop.\ \ref{L_R_plepM} is the following.

\begin{cor} \label{Hclasses_plepM} \    
Every $\cal H$-class of ${\it plep}M_{k,1}$ (or of ${\it tlep}M_{k,1}$)
has the form $H \cap {\it plep}M_{k,1}$ (respectively
$H \cap {\it tlep}M_{k,1}$), where $H$ is an $\cal H$-class of $M_{k,1}$.
 \ \ \ $\Box$
\end{cor}
Another consequence of Prop.\ \ref{L_R_plepM} is that the $\cal R$-height 
and $\cal L$-height functions that we defined for $M_{k,1}$ 
also work for ${\it plep}M_{k,1}$ and ${\it tlep}M_{k,1}$.

\medskip

\noindent Additional facts about the $\cal R$- and $\cal L$-orders of 
${\it tlep}M_{k,1}$ and ${\it plep}M_{k,1}$:

Every $\cal R$-class of $M_{k,1}$ intersects
${\it plep}M_{k,1}$, and every non-zero $\cal R$-class of $M_{k,1}$
intersects ${\it tlep}M_{k,1}$.
In particular, for any $\varphi \in M_{k,1}$ with table
$\varphi: P \to Q$ we have
 \ $\varphi \equiv_{\cal R} {\sf id}_Q \in {\it plep}M_{k,1}$.
To find an idempotent of ${\it tlep}M_{k,1}$ in every non-zero
$\cal R$-class of $M_{k,1}$ we can just extend ${\sf id}_Q$ to a total
function (by taking a complementary prefix code of $Q$ and mapping it
to any element of $Q$).

Not every $\cal L$-class of $M_{k,1}$ contains
an element of ${\it plep}M_{k,1}$.
For example, if some class of ${\sf part}(\varphi)$ contains words of
different lengths then the $\cal L$-class of $\varphi$ does not intersect
${\it plep}M_{k,1}$.

\begin{pro} \label{J_simple_plepM} \
The monoid ${\it plep}M_{k,1}$ is {\em $0$-$\cal J$-simple} (i.e., it 
consists of {\bf 0} and one non-zero $\cal J$-class), and it is 
{\em congruence-simple} (i.e., there are only two congruences in 
${\it plep}M_{k,1}$, the equality relation, and the one-class congruence).
The monoid ${\it tlep}M_{k,1}$ is $\cal J$-simple. 
\end{pro}
{\bf Proof.} \ It was proved in \cite{BiThomMon} (Prop.\ 2.2) that 
$M_{k,1}$ is 0-$\cal J$-simple. Congruence-simplicity of $M_{k,1}$ was 
proved (incompletely) in Theorem 2.3 in \cite{BiThomMon}; a complete proof 
appears in the proof of Prop.\ \ref{congrSimpleM} in the Appendix of the 
present paper. 

For 0-$\cal J$-simplicity and congruence-simplicity of ${\it plep}M_{k,1}$ 
we observe that the proofs for $M_{k,1}$ also apply for ${\it plep}M_{k,1}$ 
since the multipliers used in those proofs belong to ${\it plep}M_{k,1}$. 

Similarly, the proof of $\cal J$-simplicity of ${\it tot}M_{k,1}$ (in 
Prop.\ 2.2 in \cite{BiThomMon}) also works for ${\it tlep}M_{k,1}$,  
 \ \ \ $\Box$

\medskip

\noindent {\bf Question:} \ Are ${\it tot}M_{k,1}$ and ${\it tlep}M_{k,1}$
congruence-simple for all (or some) $k \geq 2$ ?

%%%

\subsection{ The $\cal D$-relation of ${\it tlep}M_{k,1}$ and
${\it plep}M_{k,1}$ }

This subsection gives another unexpected application of the Bernoulli 
measure $\mu$, namely a simple characterization of the $\cal D$-relation of
${\it plep}M_{k,1}$ and of ${\it tlep}M_{k,1}$.
Recall that for a $k$-ary rational number $a / k^n$ with $a$ not divisible
by $k$ we denote the numerator $a$ by ${\sf num}(r)$.

\begin{thm} {\bf ($\cal D$-relation of ${\it plep}M_{k,1}$ and 
${\it tlep}M_{k,1}$).}
 \label{D_in_plepM} \
For any non-zero $\varphi_1, \varphi_2 \in {\it plep}M_{k,1}$ , 

\smallskip

 \ \ \ $\varphi_1 \equiv_{{\cal D}({\it plepM})} \varphi_2$ \ \ \ iff
 \ \ \ ${\sf num}\big(\mu({\sf imC}(\varphi_1))\big)$  $ = $ 
${\sf num}\big(\mu({\sf imC}(\varphi_2))\big)$ .

\medskip

\noindent The same holds for ${\it tlep}M_{k,1}$.
\end{thm}
{\bf Proof.} [$\Rightarrow$] 
 \ If $\varphi_1 \equiv_{\cal R} \varphi_2$ then (by the 
characterization of the $\cal R$-order),
${\sf Im}(\varphi_1) =_{\sf ess} {\sf Im}(\varphi_2)$, hence
$\mu({\sf imC}(\varphi_1)) = \mu({\sf imC}(\varphi_2))$ (by 
Prop.\ \ref{ess_equal_same_mu} and Prop.\  \ref{L_R_plepM}). 
Thus ${\sf num}\big(\mu({\sf imC}(\varphi_1))\big) = $ 
${\sf num}\big(\mu({\sf imC}(\varphi_2))\big)$. 

Suppose $\varphi_1 \equiv_{\cal L} \varphi_2$. Then after essential
restrictions (if necessary), and since 
$\varphi_1, \varphi_2 \in {\it plep}M_{k,1}$, we have by the 
characterization of the $\cal L$-order:
 ${\sf domC}(\varphi_1) = {\sf domC}(\varphi_2) \subseteq A^m$ \ (for some
$m>0$), and \ ${\sf part}(\varphi_1) = {\sf part}(\varphi_2)$.
Hence, $|{\sf imC}(\varphi_1)| = |{\sf imC}(\varphi_2)|$
$ = |{\sf part}_{\sf domC}(\varphi_1)| = |{\sf part}_{\sf domC}(\varphi_2)|$, 
where ${\sf part}_{\sf domC}(\varphi_i)$ denotes the restriction of 
${\sf part}(\varphi_i)$ to ${\sf domC}(\varphi_1) = {\sf domC}(\varphi_2)$, 
and $|{\sf part}_{\sf domC}(\varphi_i)|$ denotes the number of classes of 
the partition on ${\sf domC}(\varphi_1) = {\sf domC}(\varphi_2)$.
Also, $\varphi_1, \varphi_2 \in {\it plep}M_{k,1}$ implies that
${\sf imC}(\varphi_1) \subseteq A^{n_1}$ and
${\sf imC}(\varphi_2) \subseteq A^{n_2}$, for some $n_1, n_2 >0$.
It follows that $\mu({\sf imC}(\varphi_1))$ and $\mu({\sf imC}(\varphi_2))$
are of the form

\smallskip

 \ \ \ \ \
$\mu({\sf imC}(\varphi_1)) \ = \ |{\sf imC}(\varphi_1)| \times k^{-n_1}$,
 \ \ and
 \ \ $\mu({\sf imC}(\varphi_2)) \ = \ |{\sf imC}(\varphi_2)| \times k^{-n_2}$.

\smallskip

\noindent Hence, since $|{\sf imC}(\varphi_1)| = |{\sf imC}(\varphi_2)|$
we have

\smallskip

 \ \ \ \ \ $\mu({\sf imC}(\varphi_1)) \times k^{n_2} \ = \ $
 $\mu({\sf imC}(\varphi_2)) \times k^{n_1}$.

\smallskip

\noindent  After removing powers of $k$, we obtain the $k$-reduced
numerators, hence 

\smallskip

 \ \ \ \ \    
 ${\sf num}\big(\mu({\sf imC}(\varphi_1))\big) \ = \ $ 
${\sf num}\big(\mu({\sf imC}(\varphi_2))\big)$.

\smallskip

\noindent We proved that both $\equiv_{\cal R}$ and $\equiv_{\cal L}$ 
preserve ${\sf num}\big(\mu({\sf imC}(\varphi))\big)$, hence 
$\equiv_{\cal D}$ preserves ${\sf num}\big(\mu({\sf imC}(\varphi))\big)$. 
The reasoning works in the same way for ${\it tlep}M_{k,1}$.

\smallskip

\noindent [$\Leftarrow$] 
 \ Let $\varphi_1, \varphi_2 \in {\it plep}M_{k,1}$ be represented by 
maps $\varphi_1: P_1 \to Q_1$ and $\varphi_2: P_2 \to Q_2$,
where $P_1, Q_1, P_2, Q_2$ are finite prefix codes with
$Q_1 \subseteq A^{n_1}$ and $Q_2 \subseteq A^{n_2}$, and
${\sf num}(\mu(Q_1)) = {\sf num}(\mu(Q_2))$.
By Lemma 4.1 in \cite{BiRL},
$\varphi_1 \equiv_{\cal R} {\sf id}_{Q_1}$ and
$\varphi_2 \equiv_{\cal R} {\sf id}_{Q_2}$, so we only need to prove that
 \ ${\sf id}_{Q_1} \equiv_{{\cal D}{\it plepM}} {\sf id}_{Q_2}$.

We have \ $\mu(Q_1) = |Q_1| \times k^{-n_1}$ \ and
 \ $\mu(Q_2) = |Q_2| \times k^{-n_2}$.
Moreover, the assumption is that
 \ $\mu(Q_1) = N \times k^{-j_1}$ \ and \ $\mu(Q_2) = N \times k^{-j_2}$, 
for a common numerator $N > 0$ and some $j_1, j_2 \geq 0$, such that $N$ 
is not divisible by $k$. Hence, \ $|Q_1| = N \times k^{i_1}$ \ and 
 \ $|Q_2| = N \times k^{i_2}$ \ for some $i_1, i_2 \geq 0$.

Suppose that, for example, $i_1 \geq i_2$. We can essentially restrict
${\sf id}_{Q_2}$ to ${\sf id}_{Q'_2}$ where $Q'_2 = Q_2 A^{i_1 - i_2}$
$ \subseteq A^{n_2+i_1 - i_2}$. 
Now we have

\smallskip

$|Q_2'| = |Q_2| \times k^{i_1 - i_2} = $
$N \times k^{i_2} \times k^{i_1 - i_2} \ = \ |Q_1|$ .

\smallskip

\noindent So there exists a bijection $\beta: Q_1 \to Q_2'$.
Since all words in $Q_1$ have the same length, and all words in $Q_2'$ have 
the same length, we have $\beta \in {\it plep}M_{k,1}$. Of course,
${\sf id}_{Q_2}$ and ${\sf id}_{Q_2'}$ represent the same element of 
${\it plep}M_{k,1}$.  Now we have
 \ $\beta \circ {\sf id}_{Q_1}(.) \circ \beta^{-1} = {\sf id}_{Q_2'}(.)$.
Hence, ${\sf id}_{Q_1} \equiv_{{\cal D}{\it plepM}} {\sf id}_{Q_2'}$, since
${\sf id}_{Q_1} \equiv_{\cal L} \beta \circ {\sf id}_{Q_1} $
$ \equiv_{\cal R} \beta \circ {\sf id}_{Q_1} \circ \beta^{-1}$
($= {\sf id}_{Q_2'}$). 
So, \ ${\sf id}_{Q_1} \equiv_{\cal D} {\sf id}_{Q_2'}$ in 
${\it plep}M_{k,1}$.

\medskip

In case $\varphi_1, \varphi_2 \in {\it tlep}M_{k,1}$ the same reasoning 
works, except that we replace ${\sf id}_{Q_1}$ and ${\sf id}_{Q'_2}$ 
(which are not total) by $\eta_{Q_1}: A^{n_1} \to Q_1$ and 
$\eta_{Q'_2}: A^{n_2 + i_1 - i_2} \to Q'_2$, defined as follows: 
For $q_1 \in Q_1$ we let $\eta_{Q_1}(q_1) = q_1$, and for 
$x \in A^{n_1} - Q_1$ we let $\eta_{Q_1}(x) = q_{0,1}$ (where $q_{0,1}$ is 
a fixed element, chosen arbitrarily in $Q_1$).
Note that the definition of $\eta_{Q_1}$ depends on $Q_1$, $q_{0,1}$, and
$A^{n_1} - Q_1$ (however, $n_1$ is determined by $Q_1$ since 
$Q_1 \subseteq A^{n_1}$, so $A^{n_1} - Q_1$ is determined by $Q_1$).
Similarly, for $q \in Q'_2$, $\eta_{Q'_2}(q) = q$, and for 
$x \in A^{n_2 + i_1 - i_2} - Q'_2$, $\eta_{Q'_2}(x) = q'_{0,2}$ (where
$q'_{0,2}$ is a fixed element, chosen arbitrarily in $Q'_2$). 
Then $\varphi_1 \equiv_{\cal R} \eta_{Q_1}$ and 
$\varphi_2 \equiv_{\cal R} \eta_{Q'_2}$.

As above, let $\beta: Q_1 \to Q'_2$ be a bijection; we assume in addition
that $\beta(q_{0,1}) = q'_{0,2}$. We define $B: A^{n_1} \to Q'_2$ by 
$B(q_1) = \beta(q_1)$ for all $q_1 \in Q_1$; and $B(q) = q'_{0,2}$ when
$q \in A^{n_1} - Q_1$. 
We define $B': A^{n_2 + i_1 - i_2} \to Q_1$ by $B'(q_2) = \beta^{-1}(q_2)$
for $q'_2 \in Q'_2$; and $B'(q) = q_{0,1}$ ($= \beta^{-1}(q'_{0,2})$) when
$q \in A^{n_2 + i_1 - i_2} - Q'_2$. Obviously, $\eta_{Q_1}$, $\eta_{Q'_2}$,
$B$, $B' \in {\it tlep}M_{k,1}$.  It is then straightforward to check that
$B' \circ B(.) = \eta_{Q_1}(.)$, and $B \circ B'(.) = \eta_{Q'_2}(.)$. 
Hence,

\smallskip

$B' \circ B \circ \eta_{Q_1} = \eta_{Q_1} = \eta_{Q_1} \circ B' \circ B$,  
 \ \ \ \ \    
$B \circ B' \circ \eta_{Q'_2} = \eta_{Q'_2} = \eta_{Q'_2} \circ B \circ B'$,

\smallskip

$B \circ \eta_{Q_1} \circ B' = \eta_{Q'_2}$, \ \ \ \ and \ \ \ \         
$B' \circ \eta_{Q'_2} \circ B = \eta_{Q_1}$.

\smallskip

\noindent Hence,
 \ $\eta_{Q_1} \equiv_{\cal L} B \circ \eta_{Q_1} \equiv_{\cal R}  $
$B \circ \eta_{Q_1} \circ B'$ \ ($ = \eta_{Q'_2}$).
 \ So, \ $\eta_{Q_1} \equiv_{\cal D} \eta_{Q'_2}$ in ${\it tlep}M_{k,1}$.
 \ \ \ \ \ $\Box$

\begin{pro} \label{Dclassforanindex} \
For any positive integer $i$ not divisible by $k$ there exists 
$\varphi \in {\it tlep}M_{k,1}$ such that 
 \ $i \ = \ {\sf num}\big(\mu({\sf imC}(\varphi))\big)$.
\end{pro}
{\bf Proof.} For any $i>0$ there exists a fixed-length prefix code 
$Q \subset A^n$ (for some $n > \log_k i$), with $|Q| = i$. So we have 
 \ $\mu(Q) = i/k^n$. Hence, if $i$ is not divisible by $k$ and if we take 
$\varphi = {\sf id}_Q$ we obtain the result for ${\it plep}M_{k,1}$. 
To get the result for ${\it tlep}M_{k,1}$ we extend ${\sf id}_Q$ to a total
function (by taking a complementary prefix code of $Q$ and mapping it
to any element of $Q$).
 \ \ \ $\Box$

\medskip

Theorem \ref{D_in_plepM} and Prop.\ \ref{Dclassforanindex}
give a one-to-one correspondence between the non-zero $\cal D$-classes of 
${\it plep}M_{k,1}$ (and of ${\it tlep}M_{k,1}$)
and the positive integers that are not divisible by $k$. 

So the $\cal D$-relation of ${\it plep}M_{k,1}$ is {\em not} 
induced by the $\cal D$-relation of $M_{k,1}$ (since $M_{k,1}$ has only
$k-1$ non-zero $\cal D$-classes). 
In other words (since the $\cal R$- and $\cal L$-orders of 
${\it plep}M_{k,1}$ are induced by $M_{k,1}$), there are 
$\psi, \varphi \in {\it plep}M_{k,1}$ such that
 \ ${\cal R}_{M_{k,1}}(\psi) \cap {\cal L}_{M_{k,1}}(\varphi)$ $ \neq $
$\varnothing$, but \ ${\cal R}_{{\it plep}M_{k,1}}(\psi)$ $ \cap $
${\cal L}_{{\it plep}M_{k,1}}(\varphi) = $ $\varnothing$.
(Notation: \ ${\cal R}_M(x)$ and ${\cal L}_M(x)$ denote the $\cal R$-
respectively $\cal L$-class of $x$ in a monoid $M$.)

It is also interesting a look at an example.
When $k=2$ and $A = \{a,b\}$, $M_{2,1}$ has just one non-zero $\cal D$-class,
so in $M_{2,1}$ we have \ ${\bf 1} \equiv_{\cal D} {\sf id}_{\{aa,b\}}$.
Obviously, ${\bf 1}$ and ${\sf id}_{\{aa,b\}}$ belong to ${\it plep}M_{k,1}$.
We have $\mu({\sf imC}({\bf 1})) = 1$, so the $k$-ary numerator is 1. On the
other hand, 
$\mu({\sf imC}({\sf id}_{\{aa,b\}})) = \mu(\{aa,b\})  = \frac{3}{4}$, so the
$k$-ary numerator is 3. Hence, 
${\bf 1} \not\equiv_{\cal D} {\sf id}_{\{aa,b\}}$ in ${\it plep}M_{k,1}$.

\begin{defn} \label{indexDclass_plepM} 
 \ For $\varphi \in {\it plep}M_{k,1}$ (or  
${\it tlep}M_{k,1}$) with $\varphi \neq {\bf 0}$, the positive integer 
${\sf num}\big(\mu({\sf imC}(\varphi))\big)$ is called {\em the index} of 
the $\cal D$-class of $\varphi$ in ${\it plep}M_{k,1}$ 
(or ${\it tlep}M_{k,1}$).
\end{defn}
The indices range over all the positive integers that are not divisible by 
$k$. Moreover, as we saw, the index determines one non-zero $\cal D$-class 
of ${\it plep}M_{k,1}$ (or ${\it tlep}M_{k,1}$) uniquely, and vice versa. 

\bigskip

Although the characterization of $\equiv_{\cal D}$ of ${\it plep}M_{k,1}$ 
in Theorem \ref{D_in_plepM} is simple, it is hard to picture what the 
number ${\sf num}(\mu({\sf imC}(\varphi)))$ means.
The following gives perhaps a better insight.

\begin{pro} \label{numMu_card} \
For any non-zero $\varphi_1, \varphi_2 \in {\it plep}M_{k,1}$ the following
are equivalent:

\smallskip

\noindent {\bf (1)}
 \ \ \ ${\sf num}\big(\mu({\sf imC}(\varphi_1))\big) \ = \ $
${\sf num}\big(\mu({\sf imC}(\varphi_2))\big)$

\smallskip

\noindent {\bf (2)}
 \ \ \ there are essential class-wise restrictions $\Phi_1$, $\Phi_2$ of
$\varphi_1$, respectively $\varphi_2$, such that for some $n \geq 1$,

\smallskip

 \ \ \  \ \ \ ${\sf imC}(\Phi_1) \subseteq A^n$,
 \ \ \ ${\sf imC}(\Phi_2) \subseteq A^n$, \ \ and
 \ \ \ $|{\sf imC}(\Phi_1)| \ = \ |{\sf imC}(\Phi_2)|$.
\end{pro}
{\bf Proof.} We prove [(1) $\Rightarrow$ (2)] (the converse is obvious).
Let $\varphi_1, \varphi_2 \in {\it plep}M_{k,1}$ be represented by
maps $\varphi_1: P_1 \to Q_1$ and $\varphi_2: P_2 \to Q_2$,
where $P_1, Q_1, P_2, Q_2$ are finite prefix codes with
$Q_1 \subseteq A^{n_1}$ and $Q_2 \subseteq A^{n_2}$, and
${\sf num}(\mu(Q_1)) = {\sf num}(\mu(Q_2))$. Moreover, by assumption,
 \ $\mu(Q_1) = N \times k^{-j_1}$ \ and \ $\mu(Q_2) = N \times k^{-j_2}$,
for a common numerator $N > 0$ and some $j_1, j_2 \geq 0$, such that $N$
is not divisible by $k$. Hence, \ $|Q_1| = N \times k^{i_1}$ \ and
 \ $|Q_2| = N \times k^{i_2}$ \ for some $i_1, i_2 \geq 0$.

Suppose that, for example, $i_1 \geq i_2$. We can essentially restrict
$\varphi_2: P_2 \to Q_2$ to $\varphi'_2: P'_2 \to Q'_2$ where
$Q'_2 = Q_2 A^{i_1 - i_2} \subseteq A^{n_2+i_1 - i_2}$.
Now we have
 \ \ $|Q_2'| = |Q_2| \times k^{i_1 - i_2} = $
$N \times k^{i_2} \times k^{i_1 - i_2} = |Q_1|$.
 \ \ \ $\Box$

%%%

\subsection{ The maximal subgroups of ${\it plep}M_{k,1}$ and
${\it tlep}M_{k,1}$ }

It is well known in semigroup theory that all the maximal subgroups
in the same $\cal D$-class are isomorphic, and
that the maximal subgroups are exactly the $\cal H$-classes that contain 
an idempotent.
So, to find all the maximal subgroups (up to isomorphism) we only need to 
find one idempotent (and its $\cal H$-class) in every $\cal D$-class. 

In \cite{BiFact} we defined the subgroup ${\it lp}G_{k,1}$ of 
length-preserving elements of the Thompson-Higman group $G_{k,1}$; one
motivation for studying ${\it lp}G_{k,1}$ is that $G_{k,1}$ is a Zappa-Szep 
product of ${\it lp}G_{k,1}$ and $F_{k,1}$ (proved in \cite{BiFact}).
More generally, for the Higman group $G_{k,m}$ we can define the subgroup
$$ {\it lp}G_{k,m} \ = \ \{ \varphi \in G_{k,m} : \    
   (\forall x \in {\sf Dom}(\varphi)) \ |\varphi(x)| = |x| \}.$$
Note that when \ $n \equiv m$ {\sf mod} $k-1$ \ then
${\it lp}G_{k,n} \simeq {\it lp}G_{k,m}$. This is proved in the same way
as Prop.\ 3.1 in \cite{BiJD} (which shows that $n \equiv m$ {\sf mod} $k-1$ 
implies $M_{k,n} \simeq M_{k,m}$).

By Prop.\ 3.1 and Theorem 2.1 in \cite{BiJD} (and their proofs) we have:

\begin{lem} \label{G_id_P} \    
Let $P \subset A^*$ be any finite prefix code such that 
$m \equiv |P|$ {\sf mod} $k-1$, let ${\sf id}_P$ be the partial identity 
on $PA^*$, and let

\smallskip
 
$G({\sf id}_P) \ = \ \{ \varphi \in M_{k,1} \ : \ $
${\sf Dom}(\varphi) =_{\sf ess} {\sf Im}(\varphi) =_{\sf ess} PA^*$,
          and $\varphi$ is {\em injective}$\}$.

\smallskip

\noindent Then $G({\sf id}_P)$ is isomorphic to the Higman group $G_{k,m}$. 
 \ \ \ $\Box$
\end{lem} 
Recall that in this paper, the integers modulo $k-1$ are taken in the range
$\{1, \ldots, k-1\}$.

\begin{pro} \label{maxSubGroups}.

\noindent {\bf (1)} \ The {\em group of units} of both ${\it plep}M_{k,1}$ 
and ${\it tlep}M_{k,1}$ is ${\it lp}G_{k,1}$. 

\medskip

\noindent {\bf (2)} \ The {\em maximal subgroups} of ${\it plep}M_{k,1}$ 
(and of ${\it tlep}M_{k,1}$) are isomorphic to the groups
${\it lp}G_{k,m}$ (for $1 \leq m \leq k-1$). More precisely, for any 
positive integer $i$ not divisible by $k$, all the maximal subgroups of the 
$\cal D$-class with index $i$ are isomorphic to
 \ ${\it lp}G_{k, \, i \, {\sf mod} \, k-1}$.
\end{pro}
{\bf Proof.} \ (1) By Corollary \ref{Hclasses_plepM}, every $\cal H$-class 
of ${\it plep}M_{k,1}$ is of the form $H \cap {\it plep}M_{k,1}$, where $H$ 
is any $\cal H$-class in $M_{k,1}$. The group of units of 
${\it plep}M_{k,1}$ is the $\cal H$-class of {\bf 1} in ${\it plep}M_{k,1}$, 
and the $\cal H$-class of {\bf 1} in  $M_{k,1}$ is $G_{k,1}$ (by Prop.\ 2.1 
in \cite{BiThomMon}). Hence, the group of units of ${\it plep}M_{k,1}$ is 
$G_{k,1} \cap {\it plep}M_{k,1} = G_{k,1} \cap {\it tlep}M_{k,1} = $
${\it lp}G_{k,1}$. 

Since the $\cal L$-class of {\bf 1} in $M_{k,1}$ contains only elements with 
domain essentially equal to $A^*$, the group of units of ${\it plep}M_{k,1}$ 
is in ${\it tlep}M_{k,1}$. Hence, the groups of units of ${\it tlep}M_{k,1}$ 
is equal to the group of units of ${\it plep}M_{k,1}$.
This proves part (1) of the Theorem.

\smallskip

\noindent Proof of (2) for ${\it plep}M_{k,1} :$ 

Let $D_i$ be the $\cal D$-class of ${\it plep}M_{k,1}$ with index $i$.
We choose any $n > 0$ such that $i < k^n$, and any prefix code 
$Q \subset A^n$ such that $|Q| = i$. Then the partial identity ${\sf id}_Q$
is an idempotent in $D_i$. The $\cal H$-class of ${\sf id}_Q$ in $M_{k,1}$ 
consists of the elements $\varphi \in M_{k,1}$ such that 
$\varphi \equiv_{\cal R} {\sf id}_Q$ (i.e., 
${\sf Im}(\varphi) =_{\sf ess} QA^*$), and 
$\varphi \equiv_{\cal L} {\sf id}_Q$ (i.e., 
${\sf Dom}(\varphi) =_{\sf ess} QA^*$, and $\varphi$ is injective). 
Hence, the $\cal H$-class of ${\sf id}_Q$ in $M_{k,1}$ is

\medskip

 \ $G({\sf id}_Q) \ = \ \{ \varphi \in M_{k,1} : \ $
 ${\sf Dom}(\varphi) =_{\sf ess} {\sf Im}(\varphi) =_{\sf ess} QA^*$, 
          and $\varphi$ is injective$\}$ , 

\medskip

\noindent and by Lemma \ref{G_id_P} this is a group isomorphic to $G_{k,m}$.
By Corollary \ref{Hclasses_plepM}, the $\cal H$-class of ${\sf id}_Q$ in 
${\it plep}M_{k,1}$ is $G({\sf id}_Q) \cap {\it plep}M_{k,1}$, hence it is 
isomorphic to ${\it lp}G_{k,m}$. 
This proves (2) for ${\it plep}M_{k,1}$.

\medskip

\noindent Proof of (2) for ${\it tlep}M_{k,1} :$ 

Let $Q \subset A^n$ be as in the proof  of (2) for ${\it plep}M_{k,1}$, 
and let $q_0$ be a fixed element, arbitrarily chosen in $Q$. 
Consider the idempotent $\eta_{Q,q_0} \in {\it tlep}M_{k,1}$ defined 
by $\eta_{Q,q_0}(q) = q$ for all $q \in Q$, and $\eta_{Q,q_0}(x) = q_0$ 
for all $x \in A^n -Q$.
The $\cal H$-class of $\eta_{Q,q_0}$ in ${\it tlep}M_{k,1}$ consists of the
elements $\varphi \in M_{k,1}$ such that we have
$\varphi \equiv_{\cal R} \eta_{Q,q_0}$ (i.e., 
${\sf Im}(\varphi) =_{\sf ess} QA^*$), and we have
$\varphi \equiv_{\cal L} \eta_{Q,q_0}$ (i.e., 
${\sf Dom}(\varphi) =_{\sf ess} A^*$, and ${\sf part}(\varphi)$ is
essentially equivalent to the partition
 $\{ \{q\} : q \in Q\} \cup \{ A^n -Q\} \}$ of $A^n$). 
So the $\cal H$-class of $\eta_Q$ in $M_{k,1}$ is 

\medskip

 $G(\eta_{Q,q_0}) \ = \ $
 \ $\big\{ \varphi \in M_{k,1} : \ \ {\sf Dom}(\varphi) =_{\sf ess} A^*$,
 \ \ ${\sf Im}(\varphi) =_{\sf ess} QA^*$, 

\hspace{1.75in}
$\varphi$ is injective on $QA^* \cap {\sf Dom}(\varphi)$, 

\hspace{1.75in} and \ $\varphi(A^n - Q) = \{ \varphi(q_0)\} \, \big\}$. 
 
\medskip

\noindent 
Then $G(\eta_{Q,q_0})$ is a maximal subgroup of $M_{k,1}$ with
identity $\eta_{Q,q_0}$ (since the $\cal H$-class of an idempotent is a 
maximal subgroup).

We saw that for $\varphi \in G(\eta_{Q,q_0})$, the restriction 
 \ $\varphi_Q:$
$QA^* \cap {\sf Dom}(\varphi) \ \to \ QA^* \cap {\sf Im}(\varphi)$ 
 \ is injective. Hence, the inverse of $\varphi \in G(\eta_{Q,q_0})$ is 
 \ $\varphi' = $ $\varphi_Q^{-1} \circ \eta_{Q,q_0}(.) \in G(\eta_{Q,q_0})$;
indeed, one easily verifies that
\ $\varphi' \varphi = \varphi \varphi' = \eta_{Q,q_0}$.

\bigskip

\noindent {\sf Claim:} 
 \ $G(\eta_{Q,q_0})$ is isomorphic to $G({\sf id}_Q)$, and
$G(\eta_{Q,q_0}) \cap {\it tlep}M_{k,1}$
is isomorphic to $G({\sf id}_Q) \cap {\it plep}M_{k,1}$.

\smallskip

\noindent
An isomorphism can be defined by \ $\iota: \varphi \in G(\eta_{Q,q_0}) $
$ \ \longmapsto \ \varphi_Q \in G({\sf id}_Q)$, where $\varphi_Q$ is the 
the restriction \ $QA^* \cap {\sf Dom}(\varphi) \ \to \ QA^* \cap {\sf
Im}(\varphi)$ \ of $\varphi$ as above.  
Bijectiveness and the homomorphism property of the map $\iota$ follow 
easily from the definition of $G(\eta_{Q,q_0})$ and $G({\sf id}_Q)$.

The same isomorphism $\iota$ shows that 
$G(\eta_{Q,q_0}) \cap {\it tlep}M_{k,1}$ is isomorphic to 
$G({\sf id}_Q) \cap {\it plep}M_{k,1}$, i.e., to ${\it lp}G_{k,m}$. 
This proves the Claim.

\medskip

By Corollary \ref{Hclasses_plepM}, $G(\eta_{Q,q_0}) \cap {\it tlep}M_{k,1}$
is the $\cal H$-class of $\eta_{Q,q_0}$ in ${\it tlep}M_{k,1}$.
This proves (2) for ${\it tlep}M_{k,1}$.
 \ \ \ \ \ $\Box$

%%%%%%%%%%%%%%%%%%%%%%%%%%%%%%%%%%%%%%%%%%%%%%%%%%%
%% Section
%%%%%%%%%%%%%%%%%%%%%%%%%%%%%%%%%%%%%%%%%%%%%%%%%%%

\section{ Complexity of computing the Bernoulli measure }

We consider the problem of computing the numbers 
$\mu({\sf imC}(\varphi))$ ($= {\sf height}_{\cal R}(\varphi)$),
$\mu({\sf domC}(\varphi))$, and the amount of collision
${\sf coll}(\varphi)$  ($= 1 - {\sf height}_{\cal L}(\varphi)$).
Here we assume that the input $\varphi \in M_{k,1}$ is given by a word
over a generating set of $M_{k,1}$; we can consider either a finite 
generating set $\Gamma$, or a circuit-like generating set 
$\Gamma \cup \tau$.
These numbers belong to $[0,1] \cap {\mathbb Z}[\frac{1}{k}]$, and we want 
to express them in base $k$, i.e., in the form 0, or 1, or 
$0.d_1 \ldots d_n$ with $d_n \neq 0$, where 
$d_1, \ldots, d_n$ $\in \{0,1, \ldots, k-1\}$ are base-$k$ digits.
Before we compute these numbers be need some preliminary algorithms.

\subsection{ Complexity for inputs over a finite generating set }

\begin{lem} \label{Compute_mu_phiInv_Gamma} \  
The following computational problems can be solved in deterministic 
polynomial time. \\  
{\sf Input:} \ $y \in A^*$ and $\varphi \in M_{k,1}$, the latter given by a 
word over a finite generating set $\Gamma$. \\  
{\sf Output 1:} \ The $k$-ary rational number $\mu(\varphi^{-1}(y))$
expressed in base $k$. \\   
{\sf Output 2:} \ The $k$-ary rational number $\mu(m_y)$ expressed in 
base $k$, where $m_y$ is any minimum-length element of the class 
$\varphi^{-1}(y)$ of ${\sf part}(\varphi)$.

The output in both cases will be a finite string over the alphabet 
$\{ \cdot, 0, 1, \ldots, k-1\}$, where ``$\cdot$'' is the base-$k$ dot.
\end{lem}
{\bf Proof.} \ We use the following result from Corollary 4.15 in
\cite{BiThomMon}:  
There is a deterministic algorithm which on input $(\varphi, y)$ 
constructs an acyclic DFA (deterministic finite automaton) ${\cal A}_y$ 
with a single accept state, that accepts the language 
$\varphi^{-1}(y) \subseteq A^*$.
The time complexity of this algorithm is a polynomial in 
$|y| + |\varphi|_{\Gamma}$ \ (where $|\varphi|_{\Gamma}$ denotes the 
word-length of $\varphi$ over $\Gamma$).

We saw that the ${\sf part}(\varphi)$-classes are finite, so the language 
$\varphi^{-1}(y)$ is finite; but its cardinality can grow exponentially 
with $|y| + |\varphi|_{\Gamma}$. On the other hand, since ${\cal A}_y$ 
can be constructed deterministically in polynomial time, ${\cal A}_y$ has 
only polynomially many states and edges. 
The underlying directed graph of ${\cal A}_y$ is acyclic, and it has only 
one source (namely the start state $q_0$) and one sink (namely the accept 
state $q_{\rm acc}$). By definition, a source
is a vertex without incoming edges, and a sink is a vertex without outgoing
edges; a finite acyclic directed graph always has at least one source and at 
least one sink. 

\smallskip

Let us compute {\sf Output 2} first: A breadth-first search is performed 
in the directed graph of ${\cal A}_y$, starting at the start state $q_0$ 
and ending when the accept state is found. This easily yields the length
of a shortest path from the start state to the accept state; and this 
length is $|m_y|$. {\sf Output 2} is $\mu(m_y) = 0.0^{|m_y|-1}1$ \ (i.e.,
after the dot in the base-$k$ representation there are $|m_y|-1$ digits
``0'' and one digit ``1'').

\medskip
  
To compute {\sf Output 1}, a more elaborate algorithm is needed.
For every state $q$ of ${\cal A}_y$ we will compute $\mu(L_q)$, where 
$L_q \subset A^*$ is the set of labels of all the paths in ${\cal A}_y$ 
from the start state $q_0$ to $q$. 
In other words, $L_q$ is the language that would be accepted by 
${\cal A}_y$ is $q$ (instead of $q_{\rm acc}$) were the {\it accept state}. 
We call $\mu(L_q)$ the {\it measure of the state} $q$, and denote it by
$\mu(q)$.
Since ${\cal A}_y$ is acyclic, $L_q$ is finite for every state $q$.
It follows from these definitions that $\mu(q_0) = 1$ (since 
$L_{q_0} = \{ \varepsilon\}$), and $L_{q_{\rm acc}} = \varphi^{-1}(y)$, so 
$\mu(q_{\rm acc}) = \mu(\varphi^{-1}(y))$.  Hence, $\mu(q_{\rm acc})$ is 
the desired {\sf Output 1}. 

\medskip

\noindent {\sf Claim.} \ For every state $q$ of ${\cal A}_y$ we have

\smallskip

 \ \ \  \ \ \ \ \      
 $\mu(q) \ = \ \frac{1}{k} \sum_{p \in {\sf pre}(q)} \mu(p)$ ,
\smallskip
 
\noindent 
where ${\sf pre}(q)$ is the set of parents (i.e., the direct predecessors)
of $q$ in the directed graph of ${\cal A}_y$. 

Indeed, \ $L_q \ = \ $
$ \bigcup \, \{ L_p \, a \ : \ q = \delta(p,a), \ p \in Q, \ a \in A\}$, 
where $Q$ denotes the set of states of ${\cal A}_y$ and 
$\delta: Q \times A \to Q$ is the next-state function.

\medskip

Based on the Claim, we obtain the following polynomial-time algorithm.
First, measure 1 is assigned to $q_0$, and measure 0 is assigned to all
other states.
Then iteratively, 
%in order of increasing depth in the acyclic directed graph of ${\cal A}_y$, 
the algorithm does the following:
A state $q$ is picked whose current measure is 0 and whose direct 
predecessors all have non-zero measure, and $q$ now receives measure 
 \ $\mu(q) = \frac{1}{k} \sum_{p \in {\sf pre}(q)} \mu(p)$.  
The algorithm never changes non-zero measures and ends when all states 
have received a non-zero measure.

\smallskip

At any moment, let ${\cal G}_0$ be the directed subgraph of ${\cal A}_y$ 
spanned by the vertices with measure 0. Since ${\cal G}_0$ is a subgraph of
an acyclic graph, ${\cal G}_0$ is acyclic; hence, ${\cal G}_0$ has sources
(i.e., vertices with no incoming edges). The predecessors in ${\cal A}_y$ of
these sources are outside of ${\cal G}_0$, so they have non-zero measure. 
The algorithm gives non-zero measure to one of these sources (i.e., it 
removes this source from ${\cal G}_0$). Since an acyclic directed graph 
always has at least one source, the algorithm will continue until all 
vertices of ${\cal A}_y$ have received a non-zero measure (i.e., until 
${\cal G}_0$ is empty). 

Breadth-first search could be used (with the queue always containing the
sources of ${\cal G}_0$), to make the algorithm more efficient.
 \ \ \ $\Box$

\begin{thm} \label{Compute_mu_Gamma} \  
On input $\varphi \in M_{k,1}$, given by a word over a finite
generating set $\Gamma$, the following numbers (expressed in
base $k$) can be computed by a deterministic polynomial-time algorithm: 
 \ $\mu({\sf imC}(\varphi))$ (i.e., ${\sf height}_{\cal R}(\varphi)$), 
 \ $\mu({\sf domC}(\varphi))$, and 
 \ ${\sf coll}(\varphi)$ (i.e., $1 - {\sf height}_{\cal L}(\varphi)$).
\end{thm}
{\bf Proof.} (1) For $\mu({\sf imC}(\varphi))$ we use the following 
result from Corollary 4.11 in \cite{BiThomMon}.
There is a deterministic polynomial-time algorithm which on input $\varphi$
(expressed over $\Gamma$), outputs ${\sf imC}(\varphi)$ explicitly as a 
list of words.

Next, for a given word $w \in A^*$ the measure $\mu(w)$ can be immediately
computed: \ $\mu(w) = 0.0^{|w|-1}1$, where $0^{|w|-1}$ denotes a 
sequence of $|w|-1$ zeros. Thus, we obtain  
 \ $\mu({\sf imC}(\varphi)) = \sum_{y \in {\sf imC}(\varphi)} \mu(y)$
 \ deterministically in polynomial-time. 

\smallskip

(2) By part (1) of this proof and from Output 1 of Lemma 
\ref{Compute_mu_phiInv_Gamma} we can compute the sequence 
$\big( \mu( \varphi^{-1}(y)) : y \in {\sf imC}(\varphi) \big)$ explicitly
in deterministic polynomial time. Hence we find 
 \ $\mu({\sf domC}(\varphi)) = $
  $\sum_{y \in {\sf imC}(\varphi)} \mu( \varphi^{-1}(y))$ \ in
polynomial time.

\smallskip

(3) For ${\sf coll}(\varphi)$ we use part (1) of this proof and 
Output 2 of Lemma \ref{Compute_mu_phiInv_Gamma} to compute
 \ ${\sf coll}(\varphi) = 1 - \sum_{y \in {\sf imC}(\varphi)} \mu(m_y)$. 
 \ \ \ $\Box$

\bigskip

For a $k$-ary rational number $0.d_1 \ldots d_n$ (with $d_n \neq 0$) 
written in base $k$, the $k$-reduced {\it numerator} is 
${\sf num}(0.d_1 \ldots d_n)  = d_1 \ldots d_n$; the latter may have
leading zeros (if $d_1 = 0$, etc.) that will be dropped.
Hence Theorem \ref{Compute_mu_Gamma} implies the following.

\begin{cor} \label{computingNum} \ If $\varphi \in M_{k,1}$ 
is given by a word over a finite generating set $\Gamma$, the following 
integers (expressed in base $k$) can be computed by a deterministic 
polynomial-time algorithm:
 \ ${\sf num}\big(\mu({\sf imC}(\varphi))\big)$, 
 \ ${\sf num}\big(\mu({\sf domC}(\varphi))\big)$, and
 \ ${\sf num}\big(\mu({\sf coll}(\varphi)\big)$.  \ \ \ $\Box$
\end{cor}

%%%%%%%%%%%%

\subsection{ New complexity classes for the Bernoulli measure and for
  counting }

Measures are similar to counting, so along the same lines as the
counting complexity classes we can define measure classes.

\begin{defn} \label{muDotC} \ 
For any complexity class $\cal C$ of decision problems we introduce the 
{\em measure class} $\mu \bullet {\cal C}$ consisting of all functions of 
the form 

\smallskip

$f_R : \ \  v \in B^* \ \longmapsto \ $
$\mu \big( \{ w \in A^* : (v,w) \in R\} \big) \ \in \ $
$[0,1] \cap {\mathbb Z}[\frac{1}{k}]$ , 

\smallskip 

\noindent where $R$ ranges over all predicates 
$R \subseteq B^* \times A^*$ (for any finite alphabets $A,B$ with $|A|=k$), 
with the following properties: 

\noindent $\bullet$ \ The predicate $R$ is {\em polynomially balanced};  
by definition, this means that there exists a polynomial $p(.)$ such that 
for all $(v,w) \in R$, \ $|w| \leq p(|v|)$.

\noindent $\bullet$ \ The membership problem of $R$ (i.e., the question, 
``given $(v,w) \in B^* \times A^*$ is $(v,w) \in R$ ?'') is in the 
complexity class $\cal C$.

\noindent $\bullet$ \ For every $v \in B^*$ the set
$(v)R = \{ w \in A^* : (v,w) \in R\}$ is a {\em finite prefix code}. 
(Finiteness of $(v)R$ already follows from polynomial balancedness of $R$.)
\end{defn}
Compare this with the well-known {\em counting class} $\# \bullet {\cal C}$ 
consisting of all functions of the form

\smallskip

$f_R : \ \  v \in B^* \ \longmapsto \ $
$| \{ w \in A^* : (v,w) \in R\} | \ \in \ {\mathbb N}$ .

\smallskip

\noindent where $R$ has the same properties as in Definition \ref{muDotC},
except that $(v)R$ is just a finite set (not required to be a prefix code).

\medskip

The counting class $\# \bullet {\sf P}$ (commonly just denoted 
$\# {\sf P}$) is called Valiant's class \cite{Valiant}, and many 
well-known problems are in that class. 
An example of a function in $\mu \bullet {\sf P}$ is given by the 
following Proposition. Here we will count every transposition 
$\tau_{i-1,i} \in \tau$ as having length $|\tau_{i-1,i}| = i$.
The alphabet $A \cup \{\cdot\}$ will be used for the base-$k$ 
representation of $k$-ary rationals, where $k = |A|$ and where the letters
$a_1,a_2, \ldots, a_k$ represent the $k$-ary digits $0,1, \ldots, k-1$.

\begin{pro} \label{measDom} \  
The following function belongs to $\mu \bullet {\sf P}$: 

\smallskip

 \ \ \ \ \ $\varphi \in M_{k,1} \ \longmapsto \ \mu({\sf domC}(\varphi))$
$ \ \in \ [0,1] \cap {\mathbb Z}[\frac{1}{k}]$,   

\smallskip

\noindent where $\varphi$ is given by a word over 
$\Gamma \cup \tau$, and $\mu({\sf domC}(\varphi))$ is expressed by a 
finite string over the alphabet $A \cup \{\cdot\}$.
\end{pro}
{\bf Proof.}  We consider the predicate 
 \ $R = \{ (\varphi, x) \in (\Gamma \cup \tau)^* \times A^*: x \in $
${\sf domC}(\varphi)\}$. By Prop.\ 5.5(1) in \cite{BiRL}, 
the membership problem of this predicate (called the domain code 
membership problem) is in {\sf P}. The proof of Prop.\ 5.5(1) in \cite{BiRL}
also shows that the predicate is polynomially balanced; in fact, for any
$x \in {\sf domC}(\varphi)$ we have
$|x| \leq c \cdot |\varphi|_{\Gamma \cup \tau}$ (for some constant $c$),
so the predicate is linearly balanced. And ${\sf domC}(\varphi)$ is of 
course a prefix code. 

In order to represent $\Gamma \cup \{\tau\}$ by a finite alphabet we will
express every transposition $\tau_{i-1,i} \in \tau$ by $t^i$. So 
$\Gamma \cup \tau$ is represented by the finite alphabet 
$\Gamma \cup \{t\}$.
 \ \ \ $\Box$

\medskip

The Bernoulli measure is closely related to counting, and we would
like to explore the connection between the measure complexity classes 
$\mu \bullet {\cal C}$ and the counting complexity classes 
$\# \bullet {\cal C}$.

 We have to overcome a syntactic obstacle, namely 
the fact that the functions in counting classes output natural integers, 
whereas the functions in measure classes output rational numbers in the 
interval $[0,1]$.
We therefore introduce output reductions that consist of moving the 
base-$k$ dot; these are a special case of polynomial-time output 
reductions. We will define them next and we will show that 
  \ $\overline{\mu \bullet {\cal C}} = \overline{\# \bullet {\cal C}}$,
where overlining indicates closure under polynomial-time dot-shift 
reduction. 

\begin{defn} \label{defDotShiftRed} \  
{\bf (1)} \ Let $f_1: B^* \to A_1^*$ and $f_2: B^* \to A_2^*$ be two 
total functions.  
A {\em polynomial-time output reduction} from $f_1$ to $f_2$ 
is a polynomial-time computable total function $\rho: A_2^* \to A_1^*$  
such that \ $f_1(.) = \rho \circ f_2(.)$.

\smallskip

\noindent {\bf (2)} \  Suppose $A_1 = A_2 = A \cup \{\cdot\}$
$ = \{a_1, \ldots, a_k, \cdot\}$, where $A \cup \{\cdot\}$ is used for
the base-$k$ representation. 
Suppose that \ ${\sf Im}(f_1) \cup {\sf Im}(f_2) \ \subseteq \ $
$A^* \, \{\cdot\} \, A^* \ \cup \ A^*$, i.e., the output strings of 
$f_1$ and $f_2$ contain at most one dot.
A {\em polynomial-time dot-shift reduction} is a polynomial-time output
reduction $\rho$ from $f_1$ to $f_2$ such that for all 
$z \in A^* \{\cdot\} A^* \cup A^*$ we have: \ $z$ and $\rho(z)$ are 
identical except possibly for occurrences of $a_1$ at the left end 
(corresponding to leading 0's), occurrences of $a_1$ at the right end
(corresponding to trailing  0's), and the position (including presence or
absence) of the dot. 
Equivalently, if $z$ and $\rho(z)$ are viewed as rationals in
in base-$k$ representation, they differ only by a multiplicative factor 
$k^n$ for some integer $n$ (positive or negative or 0).
\end{defn}

The defining property of the polynomial-time dot-shift reduction can also be
expressed as follows: If in both $z$ and $\rho(z)$ one deletes the dot and 
all occurrences of $a_1$ at the right end and the left end, the same string 
is obtained from $z$ and $\rho(z)$. 

\medskip

\noindent The {\em closure} of a set of functions $\cal F$ under 
polynomial-time dot-shift reduction is the set 

\smallskip

 \ \ \ \ \ $\overline{\cal F} \ = \ \{ \rho \circ f(.) \ : \ f \in {\cal F}$ 
       and $\rho$ is a polynomial-time dot-shift reduction$\}$. 

\medskip

The following Theorem is stated abstractly for a complexity class $\cal C$
with certain properties. But we are mainly thinking of the classes {\sf P}, 
${\sf NP}$ and ${\sf coNP}$.  Each one of these classes is closed under 
intersection with languages in {\sf P}  (i.e., $L \in {\cal C}$ and
$L_0 \in {\sf P}$ implies $L \cap L_0 \in {\cal C}$), and is closed under 
polynomial-time {\em disjunctive reduction}. 
Disjunctive polynomial-time reductions are defined in \cite{HemaspOgihara}. 

\begin{thm} \label{muBarequalsNumberBar} \  
Let $\cal C$ be any complexity class of decision problems that is closed 
under intersection with languages in {\sf P}, and closed 
under disjunctive polynomial-time reduction. Then 

\medskip

\hspace{1.5in}
$\overline{\mu \bullet {\cal C}} = \overline{\# \bullet {\cal C}}$ ,

\smallskip

\noindent
where overlining indicates closure under polynomial-time dot-shift 
reduction.
\end{thm}
{\bf Proof.} \ {\bf (1)} \ To prove
$\mu \bullet {\cal C} \ \subseteq \ \overline{\# \bullet {\cal C}}$,
consider any $f \in \mu \bullet {\cal C}$. So there is a predicate 
$R \subseteq B^* \times A^*$ such that $R \in {\cal C}$, 
$R$ is polynomially balanced, 
$(v)R  = \{ w \in A^* : (v,w) \in R\}$ is a finite prefix code (for 
every $v \in B^*$), 
and $f(v) = \mu((v)R)$ (for every $v \in B^*$). 
Let  $p(.)$ be the balancing polynomial of $R$.  From $R$ we construct a 
new predicate $R' \subseteq B^* \times A^*$ defined by

\smallskip

 \ \ \ \ \ $(v,z) \in R'$ \ \ \ iff \ \ \ $|z| = p(|v|)$,  and 
there exists a prefix $w$ of $z$ such that $(v,w) \in R$ .

\smallskip

\noindent 
Hence for all $v \in B^*$, $(v)R' \subseteq A^{p(|v|)}$; it follows that 
$(v)R'$ is a fixed-length prefix code; it follows also that $p(.)$ is a 
balancing polynomial for $R'$. 

The membership problem of $R'$ is in $\cal C$. Indeed, given $(v,z)$,
the relation $|z| = p(|v|)$ can be checked in deterministic polynomial
time. Since $z$ has only linearly many prefixes, checking whether some
prefix $x$ of $z$ satisfies $(v,w) \in R$ leads to at most $|z|$ 
membership tests in $R$; this problem belongs to $\cal C$ since $\cal C$ is 
closed under disjunctive polynomial-time reduction. 

Since $(v)R'$ is a fixed-length prefix code (of length $p(|v|)$) we have:
 \ $\mu((v)R') = |(v)R'| \cdot k^{p(|v|)}$.

The right ideals $(v)R' \cdot A^*$ and $(v)R \cdot A^*$ are essentially 
equal since every $z \in (v)R'$ has a prefix in $(v)R$ and every 
$w \in (v)R$ is the prefix of an element of $(v)R'$; the latter follows 
from the fact that the elements of $(v)R'$ have length $p(|v|)$ whereas all 
elements of $(v)R$ have length $\leq p(|v|)$.
Since both $(v)R$ and $(v)R'$ are prefix codes it follows now (by 
Prop.\ \ref{ess_equal_same_mu}) that $\mu((v)R') = \mu((v)R)$.

Hence, $\mu((v)R) = |(v)R'| \cdot k^{-p(|v|)}$, and thus $f$ is obtained 
by a polynomial-time dot-shift reduction from the function 
$v \in B^* \longmapsto |(v)R|$. The latter function belongs to 
$\# \bullet {\cal C}$.

In summary, the main idea in proof {\bf (1)} is to transform each prefix 
code $(v)R$ into a fixed-length prefix code $(v)R'$, while preserving
the measure; for fixed-length prefix codes there is a simple relation 
between cardinality and measure.

\smallskip

\noindent {\bf (2)} \ To prove
 \ $\# \bullet {\cal C} \ \subseteq \ \overline{\mu \bullet {\cal C}}$
 \ consider any $f \in \# \bullet {\cal C}$. So there is a predicate 
$R \subseteq B^* \times A^*$  such that $R \in {\cal C}$, $R$ is
polynomially balanced, and for all $v \in B^*: f(v) = |(v)R|$.
Note that here, $(v)R$ is not necessarily a prefix code.

Before constructing a new predicate from $R$ we introduce an injective
encoding homomorphism $c: A^* \to A^*$, defined for all 
$a_i \in A = \{a_1, a_2, \ldots, a_k\}$ by \ $a_i \mapsto a_ia_2$ .  
Then $|c(w)| = 2 \cdot |w|$ for all $w \in A^*$. By injectiveness,
$|c((v)R)| = |(v)R|$ for all $v \in B^*$. Now we define a new predicate
$R' \subseteq B^* \times A^*$ by

\smallskip

 \ \ \ \ \ $(v,z) \in R'$ \ \ \ iff \ \ \ $|z| = 2 \cdot p(|v|)$,  and
there exists $w$ such that $(v,w) \in R$ and $z \in c(w) \ a_1^*$ , 

\smallskip

\noindent where $p(.)$ is the balancing polynomial of $R$.
The role of $a_1^*$ is to pad $c(w)$ with trailing zeros in order to 
make $z$ have length $2 \cdot p(|v|)$.  
Then $(v)R'$ is a fixed-length prefix code with 
$(v)R' \subseteq A^{2 \cdot p(|v|)}$, and $2 \cdot p(.)$ is a 
balancing polynomial for $R'$. 
The membership problem of $R'$ is in $\cal C$, for similar reasons as 
in the proof of {\bf (1)}. 
Also, for the same reason as in {\bf (1)}, 
 \ $\mu((v)R') = |(v)R'| \cdot k^{-2 \cdot p(|v|)}$, for all $v \in B^*$.

Note that for $(v,z) \in R'$ there exists exactly one $w \in (v)R$ such that
$z = c(w) \ a_1^{2 \cdot p(|v|) - 2 \cdot |w|}$, because $c(w)$ ends with 
the letter $a_2$. Hence, \ $|(v)R'| = |(v)R|$. 

Now we have \ $|(v)R| = |(v)R'| = \mu((v)R') \cdot k^{p(|v|)}$, so $|(v)R|$ 
can be computed from $\mu((v)R')$ by a dot-shift. This yields a 
polynomial-time dot-shift reduction from the function $f$ to a problem in
$\mu \bullet {\cal C}$.   
  \ \ \ \ \ $\Box$

\medskip

\noindent
{\bf Remarks. (1)} \ We have some flexibility in the way we formulate 
the assumptions on $\cal C$ in Theorem \ref{muBarequalsNumberBar} above. 
Instead of closure under disjunctive reduction we could assume that 
$\cal C$ is closed under right-concatenation with free monoids; this means
that $R \in {\cal C}$ implies $R \, A^* \in {\cal C}$. 
Here we assume that any binary predicate 
$R \subseteq B^* \times A^*$ is represented by the language 
$\{ x\$ y \in B^* \$ A^* : (x,y) \in R\}$, where \$ is a letter that does
not belong to $A \cup B$. 

\smallskip

\noindent {\bf (2)} \ The proof of Theorem \ref{muBarequalsNumberBar} above 
shows more than what we stated: The equality 
 \ $\overline{\mu \bullet {\cal C}} = \overline{\# \bullet {\cal C}}$ 
 \ is {\it effective}, in the sense since that given a predicate $R$ that
represents a function $f$ in $\mu \bullet {\cal C}$ one easily finds a 
predicate $R'$ that represents a dot-shift of $f$ that belongs to
$\# \bullet {\cal C}$, and vice versa. 

\bigskip

\noindent As a consequence of Theorem \ref{muBarequalsNumberBar} and 
Prop.\ \ref{measDom} we have:

\begin{cor} \label{domCMeasComplete} \    The function problem \
$\varphi \in M_{k,1} \mapsto \mu({\sf domC}(\varphi))$ \ is 
 \ $\overline{\# \bullet {\sf P}}$-complete (when elements of $M_{k,1}$ 
are given by words over $\Gamma \cup \tau$).
\end{cor}
{\bf Proof.} It was proved at the end of Section 6.2 in \cite{BiJD} that
the function problem \ $C \mapsto |{\sf Dom}(C)|$ \ (where $C$
ranges over partial acyclic circuits) is $\# \bullet {\sf P}$-complete. 
Hence the problem is also $\overline{\# \bullet {\sf P}}$-complete with 
respect to polynomial-time parsimonious reductions and dot-shift 
reductions.  And it follows from Prop.\ \ref{measDom} and Theorem 
\ref{muBarequalsNumberBar} that the problem is in 
$\overline{\# \bullet {\sf P}}$ 
($ \, = \overline{\mu \bullet {\sf P}} \, $).
 \ \ \ $\Box$

%%%

\subsection{ Complexity for inputs over a circuit-like generating set
             $\Gamma \cup \tau$ }

We saw that computing $\mu({\sf domC}(\varphi))$ is
$\overline{\# \bullet {\sf P}}$-complete when $\varphi \in M_{k,1}$ is 
given by a word over $\Gamma \cup \tau$. We will show now that computing 
$\mu({\sf imC}(\varphi))$ is $\overline{\# \bullet {\sf NP}}$-complete, 
and that computing the amount of collision ${\sf coll}(\varphi)$ is
$\overline{\# \bullet {\sf coNP}}$-complete.  
It is known from \cite{TodaDiss} (see also \cite{HemaspVollmer}) that 
 \ $\# \bullet {\sf NP} \subseteq \# \bullet {\sf coNP}$.

Recall that in the definition of the word-length 
$|\varphi|_{\Gamma \cup \tau}$ we use
$|\gamma| = 1$ for $\gamma \in \Gamma$ and $|\tau_{i-1,i}| = i$.

\begin{lem} \label{lenInput} \ 
For all $\varphi \in M_{k,1}$ and all $y \in {\sf Im}(\varphi)$ there 
exists $x \in \varphi^{-1}(y)$ such that

\smallskip

$|x| \ \leq \ |y| + c_{\Gamma} \cdot |\varphi|_{\Gamma \cup \tau}$ , 

\smallskip

\noindent where 
$c_{\Gamma} = {\sf max}\{ \ell(\gamma) : \gamma \in \Gamma\}$.
\end{lem}
{\bf Proof.} \ Let $\varphi = \alpha_N \ldots \alpha_1$, where 
$\alpha_N , \ldots , \alpha_1 \in \Gamma \cup \tau$. 
When a generator $\gamma \in \Gamma$ is applied to an argument, the output 
is at most $\ell(\gamma)$ letters shorter than the argument (where 
$\ell(\gamma)$ denotes the length of the longest word in the table of 
$\gamma$). Hence, the length decrease $|x| - |\varphi(x)|$ is at most 
$c_{\Gamma} \cdot |\varphi|_{\Gamma \cup \tau}$.  Hence we have
 \ $|y| \geq |x| - c_{\Gamma} \cdot |\varphi|_{\Gamma \cup \tau}$.
 \ \ \ $\Box$

\begin{thm}  \label{imageMeasureGTau} \   
The following problem is $\overline{\# \bullet {\sf NP}}$-complete. \\  
{\sf Input:} \ $\varphi \in M_{k,1}$, given by a word over 
 $\Gamma \cup \tau$. \\   
{\sf Output:} \ The $k$-ary fraction $\mu({\sf imC}(\varphi))$, written
in base $k$. 
\end{thm}
{\bf Proof.} 
 \ We first show that this problem is $\overline{\# \bullet {\sf NP}}$-hard
(or, equivalently, $\# \bullet {\sf NP}$-hard), by 
reducing the {\it image-size problem} to it. The latter problem was proved 
to be $\# \bullet {\sf NP}$-complete in Section 6.2 of \cite{BiJD}; it is
specified as follows.  \\      
Input: \ A partial circuit $C$ with $m$ input wires and $n$ output wires;
$m, n$ are part of the input. \\   
Output: \ The integer $|{\sf Im}(C)|$, written in base $k$.

To obtain a reduction we view a partial circuit $C$ as an element (call 
it $\varphi_C$) of $M_{k,1}$, written as a word over 
$\Gamma \cup \tau$.  Indeed, any partial function between finite prefix 
codes (in this case, subsets of $A^m$ and $A^n$) is an element of 
$M_{k,1}$.  For the representation of $C$ by $\varphi_C$ we have 
 \ $|{\sf Im}(C)| = |{\sf Im}(\varphi_C) \cap A^n|$.  By knowing $n$ and 
finding $\mu({\sf imC}(\varphi_C))$ ($ = \mu({\sf Im}(\varphi_C) \cap A^n)$),
we can solve the image size problem, by computing  $|{\sf Im}(C)|$ 
$ =  k^n \cdot \mu({\sf imC}(\varphi_C))$. 
This is a polynomial-time dot-shift reduction.

\smallskip

\noindent Proof that that the problem is in 
$\overline{\# \bullet {\sf NP}} :$
 \ We consider the predicate $R \subset A^* \times (\Gamma \cup \tau)^*$ 
defined by 

\smallskip

$(y, \varphi) \in R$ \ \ iff \ \ $(\exists x \in A^*)[ \, y = \varphi(x)$
{\sf and} $|y| = c_{\Gamma} \cdot |\varphi| \, ]$, 

\smallskip

\noindent where 
$c_{\Gamma} = {\sf max}\{ \ell(\gamma) : \gamma \in \Gamma\}$.
The quantified variable $x$ in $(\exists x \in A^*)$ has polynomially
bounded length.  Indeed,
if $|y| = c_{\Gamma} \cdot |\varphi|$ and $y \in {\sf Im}(\varphi)$ then 
there exists $x \in A^*$ such that $y = \varphi(x)$ and (by Lemma
\ref{lenInput}), $|x| \leq |y| + $
$|\varphi| \cdot {\sf max}\{ \ell(\gamma) : \gamma \in \Gamma\}$.
Also, $y = \varphi(x)$ can be verified in polynomial time when $x$, $y$ and
$\varphi$ are given.
So, the membership problem of the predicate $R$ is in {\sf NP}. 

Hence the function 
 \ $\varphi \ \longmapsto \ $
$|\{ y : (y, \varphi) \in R\}|$ \ is in $\# \bullet {\sf NP}$; here, 
$\varphi$ is represented by a word over $\Gamma \cup \tau$ and the 
integer $|\{ y : (y, \varphi) \in R\}|$ is written in base $k$. 

By Theorem 4.5(2) in \cite{BiThomMon}, the length $\ell(\varphi)$ of the 
longest word in ${\sf imC}(\varphi) \cup {\sf domC}(\varphi)$ is at most 
$c_{\Gamma} \cdot |\varphi|$. Hence, if $n = c_{\Gamma} \cdot |\varphi|$ 
then the right ideal ${\sf Im}(\varphi) \cap A^nA^*$ is essential in 
${\sf Im}(\varphi)$.
Hence, ${\sf Im}(\varphi) \cap A^n$ is equal to ${\sf imC}(\Phi)$ for 
some essentially equal restriction $\Phi$ of $\varphi$.
Moreover, since ${\sf imC}(\Phi) \subseteq A^n$, ${\sf imC}(\Phi)$ is a 
prefix code. By the definition of the predicate $R$ we have 
$|\{ y : (y, \varphi) \in R\}| = |{\sf imC}(\Phi)|$;  thus, 
$|{\sf imC}(\Phi)|$ is computable in $\# \bullet {\sf NP}$.
Finally, from $|{\sf imC}(\Phi)|$ we can compute 
$\mu({\sf imC}(\varphi)) = \mu({\sf imC}(\Phi))$ in deterministic 
polynomial time, by computing 
 \ $\mu({\sf imC}(\Phi)) = k^{-n} \cdot |{\sf imC}(\Phi)|$. 
This is a polynomial-time dot-shift reduction, hence 
$\mu({\sf imC}(\Phi))$ is computable in polynomial time from 
$|{\sf imC}(\Phi)|$, which itself is computable in $\# \bullet {\sf NP}$. 
 \ \ \ $\Box$

\medskip

For the next theorem we will need a lemma.

\begin{lem} \label{allXexistsZ} \ 
Let $B(x,y)$ be any boolean formula $B(x,y)$ where $x$ and $y$ are strings 
of boolean variables with $|x| = m$ and $|y| = n$.
Suppose that there exists $y \in \{0,1\}^n$ such that 
 \ $(\forall x)[B(x,y)=0]$; i.e., $B(x,y)$ does {\em not} have property 
{\bf (2)} below.  Then there exists a boolean
formula $\beta(X,y)$ with $|X| = m+1$ and $|y| = n$, such that: 

\smallskip

\noindent {\bf (1)} \ \ $\{ y : (\forall x)[B(x,y)=1]\} \ = \ $
$\{ y : (\forall X)[\beta(X,y)=1]\}$ ;

\smallskip

\noindent {\bf (2)} \ \ the sentence
 \ $(\forall y)(\exists X)[\beta(X,y)=1]$ 
 \ is true; 

\smallskip

\noindent {\bf (3)} \ \ $\beta$ can be constructed from $B$ in
deterministic polynomial time.
\end{lem}
{\bf Proof.} We will write $x = (x_1, \ldots, x_m)$, and
$X = (x, x_{m+1})$. 
Assuming $B(x,y)$ does not have property (2) already, we define       
$\beta(x, x_{m+1}, y)$ by 

\medskip

 \ \ \ $\beta(x, 0, y) \ = \ B(x,y)$ , \ \ and

\smallskip

 \ \ \ $\beta(x, 1, y) \ = \ 1$ . 

\medskip

\noindent Then for every $y$ we have $(\forall x)[\beta(x, 1, y) = 1]$, 
hence $(\exists x)[\beta(x, x_{m+1}, y) =1]$; so property (2) holds.
The definition of $\beta(x, x_{m+1}, y)$ can also be written as

\smallskip

 \ \ \ $\beta(x, x_{m+1}, y) \ = \ x_{m+1} \ \vee \ B(x,y)$ , 

\smallskip

\noindent which provides an expression for $\beta(x, x_{m+1}, y)$ from
an expression for $B(x,y)$ in deterministic linear time. This proves 
property (3). Let us prove of property {\bf (1)}.

\smallskip

\noindent [$\subseteq$] \ \ If $y$ satisfies $(\forall x)[B(x,y)=1]$ then 
$(\forall x)[\beta(x,0,y)=1]$ and $(\forall x)[\beta(x,1,y)=1]$. Hence
$y$ satisfies \ $(\forall (x,x_{m+1}))[\beta(x,x_{m+1},y)=1]$.

\noindent [$\supseteq$] \ \ If $y$ satisfies 
$(\forall (x,x_{m+1}))[\beta(x,x_{m+1},y)=1]$ then when $x_{m+1} = 0$ we 
have $(\forall x)[\beta(x,0,y) = 1]$. Hence, since
$\beta(x, 0, y) = B(x,y)$, $y$ satisfies $(\forall x)[B(x,y) = 1]$.
 \ \ \ \ \ $\Box$

\medskip

\noindent {\bf Notation.} \ For a boolean formula $B(x,y)$ we let 

\smallskip

 \ \ \ $N_{B,1} \ = \ |\{y : (\forall x)[B(x,y)=1]\}|$ , 

\smallskip

 \ \ \ $N_{B,0} \ = \ |\{y : (\forall x)[B(x,y)=0]\}|$ .

\smallskip

\noindent Hence, $N_{B,0} = 0$ \ iff \ $B(x,y)$ satisfies 
 \ $(\forall y)(\exists x)[B(x,y)=1]$ 
(i.e., property (2) in Lemma \ref{allXexistsZ}).
 
\begin{thm}  \label{collGTau} \
The following problem is $\overline{\# \bullet {\sf coNP}}$-complete. \\
{\sf Input:} \ $\varphi \in M_{k,1}$, given by a word over
 $\Gamma \cup \tau$. \\
{\sf Output:} \ The $k$-ary fraction ${\sf coll}(\varphi)$, written
in base $k$.
\end{thm}
{\bf Proof. (1)} \ Let us first prove that the function 
$\varphi \mapsto {\sf noncoll}(\varphi)$ belongs to 
$\mu \bullet {\sf coNP}$  ($\subseteq \overline{\# \bullet {\sf coNP}}$).
Recall that ${\sf coll}(\varphi) = 1 - {\sf noncoll}(\varphi)$ and 
${\sf noncoll}(\varphi) = \sum_{i=1}^n \mu(m_i)$ where 
$n = |{\sf imC}(\varphi)|$ and $\{m_i: i = 1, \ldots, n\}$ consists of
minimum-length representative of all the ${\sf part}(\varphi)$-classes 
in ${\sf domC}(\varphi)$ (according to Def.\ \ref{collision}); this assumes 
that ${\sf imC}(\varphi)$ is a prefix code. 
By Theorem 4.12 in \cite{BiThomMon} and its proof, every word in 
${\sf domC}(\varphi)$ has length 
$\leq c_{\Gamma} \cdot |\varphi|_{\Gamma \cup \tau}$, where $c_{\Gamma}$
$ = {\sf max}\{ |z| : z \in {\sf domC}(\gamma) \cup {\sf imC}(\gamma), $
$\gamma \in \Gamma \}$; i.e., $c_{\Gamma}$ is the length of the longest 
words occurring in the tables of the elements of $\Gamma$.
In this proof we will abbreviate $c_{\Gamma}$ by $c$ and 
$|\varphi|_{\Gamma \cup \tau}$ by $|\varphi|$.

When we choose a minimum-length representatives $m$ in a class $C$ we will
choose $m$ to be the first element in the dictionary order (among the 
minimum-length elements in $C$). This does not affect the value of
${\sf noncoll}(\varphi)$, which depends only on lengths. 
The {\em dictionary order} in $A^*$ is defined as follows for all 
$x, y \in A^*$
(assuming $A$ is ordered as $a_1 < \ldots < a_k$):
 \ \ $x \leq_{\sf dict} y$ \ \ iff

\smallskip

 \ \ $\bullet$ \ $x$ is a prefix of $y$, \ or

 \ \ $\bullet$ \ there is a common prefix $p$ of $x$  
 and $y$ and there are $a_i, a_j \in A$ such that  

 \ \ \ \ \  $pa_i$ is a prefix of $x$ and $pa_j$ is a prefix of $y$, 
          and $a_i < a_j$.

\smallskip

\noindent Since $x$ has only linearly many prefixes, the property 
$x \leq_{\sf dict} y$ can be checked in deterministic polynomial time when 
$x$ and $y$ are given.

By Lemma \ref{classw_restrVSlength}(2), $\varphi$ has an essential 
restriction $\Phi$ that is such that ${\sf imC}(\Phi)$ is a prefix code 
and all the minimum-length representatives $m_i$ of the 
${\sf part}(\Phi)$-classes have (the same) length $c \cdot |\varphi|$. 
By Lemma \ref{coll_vs_essequ}, 
 \ ${\sf noncoll}(\varphi) = {\sf noncoll}(\Phi)$. 

A set of minimum-length representatives of the ${\sf part}(\Phi)$-classes
in ${\sf domC}(\Phi)$ is the same as a set of representatives 
of those ${\sf part}(\varphi)$-classes $C$ such that $C$ contains some
element(s) of length $c \cdot |\varphi|$, and $C$ contains no shorter 
elements.  Therefore,

\smallskip

${\sf noncoll}(\Phi) \ = \ $
$\mu\big(\{ m \in A^{c \, |\varphi|} \cap {\sf Dom}(\varphi) \ : \ $

\hspace{1.3in}  
$(\forall x \in A^{< c \, |\varphi|}) [\varphi(x) \neq \varphi(m)]$ 
 \ \ and 
 \ \ $(\forall x \in A^{c \, |\varphi|})$
 $[x <_{\sf dict} m \Rightarrow \varphi(x) \neq \varphi(m)] \, \} \big)$.

\smallskip

\noindent Here, 
$(\forall x \in A^{< c \, |\varphi|}) [\varphi(x) \neq \varphi(m)]$ 
 \ expresses that $m$ has minimum length in its ${\sf part}(\varphi)$-class;
and \ $(\forall x \in A^{c \, |\varphi|})$
 $[x <_{\sf dict} m \Rightarrow \varphi(x) \neq \varphi(m)]$ \ expresses that
$m$ comes first in the dictionary order in its ${\sf part}(\varphi)$-class
(which implies that we have picked only one representative in this class).

The above formula for ${\sf noncoll}(\Phi)$ implies that the function 
${\sf noncoll}(.)$ belongs to $\mu \bullet {\sf coNP}$. Indeed, the 
restriction that $m \in A^{< c \, |\varphi|}$ implies polynomial 
balancedness.
The property $[m \in {\sf Dom}(\varphi)]$ is in {\sf P} (by Prop.\ 5.5 
in \cite{BiRL}). The condition $[\varphi(x) \neq \varphi(m)]$ is in
{\sf P}, since $\varphi$ can be evaluated on two given words in 
deterministic polynomial time (by the proof of Prop.\ 5.5 in \cite{BiRL});
hence the condition 
 \ $(\forall x \in A^{< c \, |\varphi|}) [\varphi(x) \neq \varphi(m)]$ 
 \ is in {\sf coNP}. 
Similarly, $[x <_{\sf dict} m \Rightarrow \varphi(x) \neq \varphi(m)]$ 
is in {\sf P}, hence the condition 
 \ $(\forall x \in A^{c \, |\varphi|})$
 $[x <_{\sf dict} m \Rightarrow \varphi(x) \neq \varphi(m)]$ \ is in 
{\sf coNP}. 

\medskip

\noindent {\bf (2)} \ Let us prove now that the function ${\sf noncoll}(.)$ 
is $\overline{\# \bullet {\sf coNP}}$-hard, by reducing the problem
$\# \bullet \Pi_1^{\sf P}${\sf Sat} to the problem of computing 
${\sf noncoll}(.)$. The problem $\# \bullet \Pi_1^{\sf P}${\sf Sat} is 
also denoted by $\# \forall${\sf Sat} and called ``Counting 
forall-satisfiability''; it is specified as follows.
 
\smallskip

\noindent {\sf Input:} \ A boolean formula $B(x,y)$ where $x$ and $y$ are 
strings of variables with $|x| = m$, $|y| = n$ (with $m$ and $n$ part of
the input). 
 
\smallskip

\noindent {\sf Output:} \ The binary representation of the integer
 \ $|\{y \in \{0,1\}^n \, : \, (\forall x \in \{0,1\}^m)[ B(x,y) = 1]\}|$.

\smallskip

\noindent
The problem $\# \forall${\sf Sat} is $\# \bullet {\sf coNP}$-complete (see
\cite{DurandHerKol}), and remains $\# \bullet {\sf coNP}$-complete when we 
restrict to the case when $n = m$; we assume from now on that $n = m$. 
For a reduction we map any instance $B(x,y)$ of $\# \forall${\sf Sat} to the
element $\varphi_B \in M_{2,1}$, defined as follows:

\medskip

 \ \ \ $\varphi_B(0xz) \ = \ B(x,z) \cdot x$ \ \ \ for all 
      $x \in \{0,1\}^n$ and $z \in \{0,1\}^n$ ;

\smallskip

 \ \ \ $\varphi_B(1xw) \ = \ 0 x$ \ \ \ for all $x \in \{0,1\}^n$
       and $w \in \{0,1\}^{n+1}$ .

\medskip

\noindent So, \ ${\sf domC}(\varphi_B) = $
$ 0 \, \{0,1\}^{2n} \cup 1 \, \{0,1\}^{2n+1}$, and
 \ $0 \, \{0,1\}^n \subseteq {\sf imC}(\varphi_B) \subseteq \{0,1\}^{n+1}$.
More precisely, \ \ ${\sf imC}(\varphi_B) \ = \ 0 \, \{0,1\}^n $ 
$ \ \cup \ $ $1 \, \{x \in \{0,1\}^n : (\exists z)[ B(z,x) = 1]\}$.
By Lemma \ref{allXexistsZ} we can assume that for every $x \in \{0,1\}^n$ 
there exists $z$ such that $B(x,z) = 1$; by the Lemma, this does not change 
the cardinality $|\{ y : (\forall x)[ B(x,y) = 1]\}|$. 
That assumption implies that \ ${\sf imC}(\varphi_B) \ = \ \{0,1\}^{n+1}$. 
 
By the definition of $\varphi_B$, the classes of ${\sf part}(\varphi_B)$ 
are of the form $\varphi_B^{-1}(0x)$ for all $x \in \{0,1\}^n$, and
of the form $\varphi_B^{-1}(1x)$ for all $x \in \{0,1\}^n$ such that 
$(\exists z)[B(x,z) = 1]$; however, this holds for all $x \in \{0,1\}^n$,
since by Lemma \ref{allXexistsZ} we assume that for all $x \in \{0,1\}^n$ we 
have $(\exists z)[B(x,z) = 1]$. 
In a class $\varphi_B^{-1}(1x)$ all elements (and hence the shortest element)
have length $2n+1$; and there are $2^n$ such classes (as $x$ ranges over
$\{0,1\}^n$).
The shortest element in the class $\varphi_B^{-1}(0y)$ has length $2n+2$ if 
$(\forall z)[B(x,z) = 1]$; the number of such classes is 
 \ $N_{B,1} =  |\{y : (\forall x)[B(x,y)=1]\}|$. 
In a class $\varphi_B^{-1}(0y)$, the shortest element has length $2n+1$ 
if $(\exists z)[B(x,z) = 0]$; the number of such classes is 
$2^n - N_{B,1}$.
Thus for ${\sf noncoll}(\varphi_B)$ we have the formula

\medskip

${\sf noncoll}(\varphi_B) \ = \ $
$2^{-(2n+1)} \cdot 2^n \ + \ 2^{-(2n+2)} \cdot N_{B,1} $
$ \ + \ 2^{-(2n+1)} \cdot (2^n - N_{B,1})$
$ \ = \ 2^{-n} - 2^{-(2n+2)} \cdot N_{B,1}$.

\medskip

\noindent It follows that, in binary representation, $N_{B,1}$ can be 
computed in deterministic polynomial time from ${\sf noncoll}(\varphi_B)$ 
via the formula

\medskip

$N_{B,1} \ = \ 2^{n+2} - {\sf noncoll}(\varphi_B) \cdot 2^{2n+2}$ ,

\medskip

\noindent which reduces the problem $\# \bullet \Pi_1^{\sf P}${\sf Sat} 
(i.e., the computation of $N_{B,1}$) to the problem of computing 
${\sf noncoll}(\varphi_B)$.        \ \ \ \ \ $\Box$

\bigskip

\bigskip

%%%%%%%%%%%%%%%%%%%%%%%%%%%%%%%%%%%%%%%%%%%%%%%%%%%
%% Section
%%%%%%%%%%%%%%%%%%%%%%%%%%%%%%%%%%%%%%%%%%%%%%%%%%%

\section{ Appendix }

The following theorem was stated in \cite{BiThomMon} (Theorem 2.3), but 
the proof was incomplete. We give a complete proof here.

\begin{thm} \label{congrSimpleM} \ 
The monoids $M_{k,1}$ and ${\it Inv}_{k,1}$ are congruence-simple 
for all $k \geq 2$.
\end{thm}
{\bf Proof.} \ Let $\equiv$ be any congruence on $M_{k,1}$ that is not 
the equality relation. We will show that then the whole monoid is 
congruent to the empty map {\bf 0}. We will make use of 
$0$-$\cal J$-simplicity.

\smallskip

\noindent {\sf Case 0:} \  
Assume that $\Phi \equiv $ {\bf 0} for some element $\Phi \neq $
{\bf 0} of $M_{k,1}$. Then for all $\alpha, \beta \in M_{k,1}$ we have 
obviously $\alpha \, \Phi \, \beta \equiv $ {\bf 0}. Moreover, by 
$0$-$\cal J$-simplicity of $M_{k,1}$ we have \ $M_{k,1}$
$= \{ \alpha \, \Phi \, \beta : \alpha, \beta \in M_{k,1}\}$
 \ since $\Phi \neq 0$. Hence in this case all elements of $M_{k,1}$ are 
congruent to {\bf 0}. 

\smallskip

\noindent For the remainder we suppose that $\varphi \equiv \psi$ and 
$\varphi \neq \psi$, for some elements $\varphi, \psi$ of 
$M_{k,1} - \{ {\bf 0} \}$.

\smallskip

\noindent {\sf Case 1:} 
 \ ${\sf Dom}(\varphi) \neq_{\sf ess} {\sf Dom}(\psi)$.

Then there exists $x_0 \in A^*$ such that 
$x_0A^* \subseteq {\sf Dom}(\varphi)$, but 
 \ ${\sf Dom}(\psi) \cap x_0A^* = \varnothing$;  or, vice versa,
there exists $x_0 \in A^*$ such that 
$x_0A^* \subseteq {\sf Dom}(\psi)$, but
 \ ${\sf Dom}(\varphi) \cap x_0A^* = \varnothing$.  Let us assume the former.
Letting $\beta = (x_0 \mapsto x_0)$, we have
$\varphi \, \beta(.) = (x_0 \mapsto \varphi(x_0))$. 
We also have $\psi \, \beta(.) = $ {\bf 0}, since 
$x_0A^* \cap {\sf Dom}(\psi) = \varnothing$.
So, $\varphi \, \beta \equiv \psi \, \beta = {\bf 0}$, but  
$\varphi \, \beta \neq$ {\bf 0}.
Hence case 0, applied to $\Phi = \varphi \, \beta$, 
implies that the entire monoid $M_{k,1}$ is congruent to {\bf 0}.

\smallskip

\noindent {\sf Case 2.1:}
 \ ${\sf Im}(\varphi) \neq_{\sf ess} {\sf Im}(\psi)$
 \ and \ ${\sf Dom}(\varphi) =_{\sf ess} {\sf Dom}(\psi)$.

Then there exists $y_0 \in A^*$ such that 
$y_0A^* \subseteq {\sf Im}(\varphi)$, but                                  
 \ ${\sf Im}(\psi) \cap y_0A^* = \varnothing$;  or, vice versa, 
$y_0A^* \subseteq {\sf Im}(\psi)$, but                
 \ ${\sf Im}(\varphi) \cap y_0A^* = \varnothing$. 
Let us assume the former. Let $x_0 \in A^*$ be such that  
$y_0 = \varphi(x_0)$. Then 
 \ $(y_0 \mapsto y_0) \circ \varphi \circ (x_0 \mapsto x_0)$
$ \ = \ (x_0 \mapsto y_0)$. 

On the other hand, 
 \ $(y_0 \mapsto y_0) \circ \psi \circ (x_0 \mapsto x_0) \ = \ {\bf 0}$. 
Indeed, if $x_0A^* \cap {\sf Dom}(\psi) = \varnothing$ then for all 
$w \in A^* : $ \ $\psi \circ (x_0 \mapsto x_0)(x_0w) \ = \ \psi(x_0w)$
$ \ = \ $ $\varnothing$. 
And if $x_0A^* \cap {\sf Dom}(\psi) \neq  \varnothing$ then 
for those $w \in A^*$ such that $x_0w \in {\sf Dom}(\psi)$ we have
 \ $(y_0 \mapsto y_0) \circ \psi \circ (x_0 \mapsto x_0)(x_0w) \ = \ $
$(y_0 \mapsto y_0)(\psi(x_0w)) \ = \ \varnothing$, since
${\sf Im}(\psi) \cap y_0A^* = \varnothing$.  Now case 0 applies to 
${\bf 0} \neq \Phi = $
$(y_0 \mapsto y_0) \circ \varphi \circ (x_0 \mapsto x_0)$
$ \equiv {\bf 0}$; 
hence all elements of $M_{k,1}$ are congruent to {\bf 0}.

\smallskip

\noindent {\sf Case 2.2:}
  \ ${\sf Im}(\varphi) =_{\sf ess} {\sf Im}(\psi)$
 \ and \ ${\sf Dom}(\varphi) =_{\sf ess} {\sf Dom}(\psi)$.

Then (after restricting), \ ${\sf domC}(\varphi) = {\sf domC}(\psi)$, and 
 there exist $x_0 \in {\sf domC}(\varphi) = {\sf domC}(\psi)$ and
$y_0 \in {\sf imC}(\varphi)$, $y_1 \in {\sf imC}(\psi)$ such
that $\varphi(x_0) = y_0 \neq y_1 = \psi(x_0)$. We have two sub-cases.

\smallskip

\noindent {\sf Case 2.2.1:} \ $y_0$ and $y_1$ are not prefix-comparable.

Then
 \ $(y_0 \mapsto y_0) \circ \varphi \circ (x_0 \mapsto x_0)$
$ \ = \ (x_0 \mapsto y_0)$.

On the other hand, 
 \ $(y_0 \mapsto y_0) \circ \psi \circ (x_0 \mapsto x_0)(x_0 w) \ = $ 
 \ $(y_0 \mapsto y_0)(y_1 w) \ = \ \varnothing$ \ for all $w \in A^*$ \ 
(since $y_0$ and $y_1$ are not prefix-comparable).  So
 \ $(y_0 \mapsto y_0) \circ \psi \circ (x_0 \mapsto x_0) \ = \ \ {\bf 0}$.
Hence case 0 applies to \ ${\bf 0} \neq\Phi \ = \ $
$(y_0 \mapsto y_0) \circ \varphi \circ (x_0 \mapsto x_0)$
$\equiv {\bf 0}$.

\smallskip

\noindent {\sf Case 2.2.2:} \ $y_1$ is a prefix of $y_0$. 
(The case where $y_0$ is a prefix of $y_1$ is similar; and since 
$y_0 \neq y_1$, there is no other case.)

Then $y_1 = y_0 v_1$ for some $v_1 \in A^*$, and 
$y_1A^* \subsetneqq y_0A^*$, so ${\sf imC}(\psi) \cap y_0A^*$ contains 
some string $y_2$ besides $y_1$.
Indeed, the right ideal ${\sf imC}(\psi) \, A^* \cap y_0A^*$ is essential 
in $y_0A^*$ because ${\sf Im}(\varphi) =_{\sf ess} {\sf Im}(\psi)$.

So, $y_2 = y_0 v_2$ for some $v_2 \in A^*$. Hence, 
 \ $(y_2 \mapsto y_2) \circ \varphi \circ (x_0 \mapsto x_0)(x_0v_2)$
$ \ = \ $ $(y_2 \mapsto y_2)(y_0 v_2) \ = \ y_2$.

On the other hand, $y_1$ and $y_2$ are not prefix-comparable, since both 
belong to ${\sf imC}(\psi)$, which is a prefix code. Hence, 
 \ $(y_2 \mapsto y_2) \circ \psi \circ (x_0 \mapsto x_0)(x_0w) \ = \ $
 $(y_2 \mapsto y_2)(y_1 w) \ = \ \varnothing$, since
$y_2$ and $y_1$ are not prefix-comparable.
Thus, case 0 applies to \ ${\bf 0} \neq \Phi = $
$(y_2 \mapsto y_2) \circ \varphi \circ (x_0 \mapsto x_0)$
$\equiv {\bf 0}$.  

\smallskip

The same proof works for ${\it Inv}_{k,1}$ since all the multipliers used
inthe proof (of the form $(u \mapsto v)$ for some $u,v \in A^*$) belong to 
${\it Inv}_{k,1}$.
 \ \ \ $\Box$

%%%%%%%%%%%%%%%%%%%%%%%%%%%%%%%%%%%%%%%%%%%%%%%%%%%%%%%%

\bigskip

\bigskip

%%%%%%%%%%%%%%%%%%%%%%%%%%%%%%%%%%%%%%%%%%%%%%%%%%%%%%%%%%%%%%%%%%%%%%%%%%%%%%%
%% \small

%%%%%%%%%%%%%%%%%%%%%%%%%%%%%

\bigskip

\bigskip

\noindent {\bf Jean-Camille Birget} \\
Dept.\ of Computer Science \\
Rutgers University at Camden \\
Camden, NJ 08102, USA \\
{\tt birget@camden.rutgers.edu}

\end{document}